\documentclass[3p, 11pt]{elsarticle}

\makeatletter
\def\ps@pprintTitle{%
 \let\@oddhead\@empty
 \let\@evenhead\@empty
 \def\@oddfoot{}%
 \let\@evenfoot\@oddfoot}
\makeatother

\usepackage{lipsum}
\usepackage{hyperref}
\usepackage{amssymb}
\usepackage{caption,subcaption}
\usepackage{amsmath}
\usepackage{pgfplots}
\usepackage{pgf}
\usepackage{tikz} 
\usepackage{amsfonts}
\usepackage{graphicx}
\usepackage{epstopdf}
\usepackage{MnSymbol}
\usepackage{overpic}
\usepackage{float}
\usepackage{graphicx}
\usepackage[absolute,overlay]{textpos}
\usepackage{booktabs}
\usepackage{relsize}
\usepackage{mathrsfs}

\usepackage[linesnumbered, ruled, vlined]{algorithm2e}
\biboptions{sort&compress}
\newcommand{\bm}[1]{\ensuremath{\mathbf{#1}}}

\newtheorem{definition}{Definition}
\newcommand{\Vhat}{\hat{V}}

\begin{document}

\begin{frontmatter}

\title{Cross Interpolation for Solving High-Dimensional Dynamical Systems on Low-Rank Tucker  and Tensor Train Manifolds}

\author{Behzad Ghahremani$^{1}$ and Hessam Babaee$^{1*}$}

\address{Department of Mechanical Engineering and Materials Science, University of
Pittsburgh, 3700 O’Hara Street, Pittsburgh, PA, 15213, USA \\ \vspace{2mm}
* Corresponding Author, Email:h.babaee@pitt.edu \vspace{-8mm}}

\begin{abstract} 
We present a novel tensor interpolation algorithm for the time integration of nonlinear tensor differential equations (TDEs) on the tensor train and Tucker tensor low-rank manifolds, which are the building blocks of many tensor network decompositions. This paper builds upon our previous work (\emph{Donello et al., Proceedings of the Royal Society A, Vol. 479, 2023} \cite{MDOblique}) on solving nonlinear matrix differential equations on low-rank matrix manifolds using CUR decompositions. The methodology we present offers multiple advantages: (i) It delivers near-optimal computational savings both in terms of memory and floating-point operations by leveraging cross algorithms based on the discrete empirical interpolation method to strategically sample sparse entries of the time-discrete TDEs to advance the solution in low-rank form. (ii) Numerical demonstrations show that the time integration is robust in the presence of small singular values. (iii) High-order explicit Runge-Kutta time integration schemes are developed. (iv) The algorithm is easy to implement, as it requires the evaluation of the full-order model at strategically selected entries and does not use tangent space projections, whose efficient implementation is intrusive. We demonstrate the efficiency of the presented algorithm for several test cases, including a nonlinear 100-dimensional TDE for the evolution of a tensor of size $70^{100} \approx 3.2 \times 10^{184}$ and a stochastic advection-diffusion-reaction equation with a tensor of size $4.7 \times 10^9$.
\end{abstract}


\end{frontmatter}

\section{Introduction} \label{Introduction}
Tensor networks can overcome the curse of dimensionality by approximating high-dimensional tensors using only lower-dimensional tensors. The density matrix renormalization group (DMRG) is the first example of a tensor network employed to represent high-dimensional quantum states, thereby mitigating the curse of dimensionality \cite{W92}. In the past decade, there has been an explosion in the development of new tensor networks in the field of quantum information sciences. See \cite{B23} for a recent review of these developments.  In recent years, tensor networks have been utilized in many diverse domains including numerical analysis, machine learning, large-scale optimization problems, kinetics, and fluid dynamics.    See \cite{CLO16} for more details on the application of tensor networks. 

The tensor train or matrix product state \cite{TT11} and the Tucker tensor  \cite{T66} are two of the widely used tensor low-rank approximation models.  The key attribute of tensor train low-rank approximation is that it does not inherently suffer from the curse of dimensionality. More specifically, tensor train low-rank approximations reduce the total number of elements in a tensor of size $n^d$ to $\mathcal{O}(dnr^2)$, where $r$ is the compression rank. For applications where $r$ is relatively small, tensor trains mitigate the curse of dimensionality due to their linear dependence on $d$. The Tucker tensor decomposition is another well-known tensor low-rank approximation model that reduces the total number of elements in a $d$-dimensional tensor from $n^d$ to $r^d + rdn$, where $r$ represents the rank of the unfolded tensors along each mode \cite{T66}. Although the Tucker model itself cannot be used for very high-dimensional tensors due to its exponential dependence on $d$, it serves as one of the building blocks for other complex tensor networks \cite{GKT13}.

One important application of tensor networks is the reduction of computational costs in solving multi-dimensional partial differential equations (PDEs) through dynamical low-rank approximation (DLRA) \cite{KL10,CL18Time}. The DLRA formalism enables the time integration of generic tensor differential equations (TDEs) on low-rank manifolds, thereby expanding the utility of tensor low-rank approximation techniques to much broader applications beyond quantum systems.  Discretizing time-dependent, multi-dimensional PDEs in all dimensions except time results in tensor differential equations (TDEs) in the form  $dV/dt = \mathcal{F}(V),$ where 
$V \in \mathbb{R}^{n_1 \times n_2 \times \dots \times n_d}$ is the solution tensor, and  $ \mathcal{F}(V)$
is the right-hand side tensor of the same size as \(V\). For applications of DLRA in solving the Schrödinger equation and Hamilton-Jacobi-Bellman equations, see \cite{AV2017, SDTensor}. The tensor train low-rank approximation has recently been used to perform direct numerical simulation (DNS) of turbulent flows by tensorizing the resulting high-dimensional dynamical system \cite{GL22}.  Solving TDEs, even in moderate dimensions (e.g., \(d=4,5,6\)), using traditional numerical methods such as finite difference and finite element methods, encounters the issue of the curse of dimensionality  \cite{IGTensor}. To decrease the number of degrees of freedom (DOF) of \(V(t)\), DLRA constrains the solution of the TDEs to the manifold of low-rank tensors, where explicit evolution equations are obtained in the low-rank form with significantly smaller DOF.

Despite the effectiveness of DLRA in reducing the computational cost of solving multi-dimensional TDEs, several challenges remain as enumerated below:

(i) \textbf{Computational Cost:} The computational cost of solving DLRA increases when the exact rank of $\mathcal{F}(V)$ is high. This issue arises in linear TDEs featuring a large number of right-hand side terms or full-rank forcing matrices \cite{DCB22,AB24} and in TDEs with high-order polynomial nonlinearities. For TDEs that incorporate non-polynomial nonlinearities, such as exponential or fractional nonlinearity, the floating-point operation (FLOP) costs of solving the TDE on low-rank manifolds exceed those of the full-order model (FOM). This increase in cost is due to the fact that $\mathcal{F}(V)$ becomes a full-rank tensor even when $V$ is of low rank. Consequently, $\mathcal{F}(V)$ needs to be computed explicitly, requiring $\mathcal{O}(n^d)$ FLOP costs, which are similar to those of the FOM. While the DLRA equations can still be solved in a memory-efficient manner by sequentially computing entries of $\mathcal{F}(V)$, this approach results in increased wall clock time due to the sequential nature of the algorithm.

(ii) \textbf{Intrusiveness:} Even when the exact rank of $\mathcal{F}(V)$ is low, for example, for linear TDEs, the efficient implementation of DLRA evolution equations is intrusive, for example, when dealing with complex codes, such as multi-physics problems. This process involves substituting the low-rank approximation in the FOM, and projection of $\mathcal{F}(V)$ onto the tangent space of the manifold \cite{MDOblique, MNAdaptive}. However, to retain the memory and FLOP efficiencies that the DLRA offers, the implementation must be done carefully to avoid forming tensors in the ambient space or forming tensors that require significantly more degrees of freedom than what is needed to represent $V$ in the low-rank form.

(iii) \textbf{Ill-Conditioning:} The DLRA evolution equations obtained from the projection onto the tangent space of the manifold may encounter numerical stability issues. The computation of the evolution equations involves inverting auto-correlation matrices that may become ill-conditioned \cite{CL15TT, CL18Time, JHUnifying,GC22Unconventional}. This is especially problematic when the rank ($r$) of the system must increase to achieve a more accurate approximation. The inversion of these possibly ill-conditioned matrices causes stringent constraints on time step size in numerical integration and error amplification.

The aforementioned issues have been the focus of intense research in the past decade.   Projector-splitting techniques \cite{CLProjector,CL15TT}  
utilize a Lie–Trotter or Strang splitting of the tangent-space projection. These time integration schemes are robust in the presence of small or zero singular values, however, they involve a backward time step, which renders these techniques unstable for dissipative dynamical systems. The robust basis update \& Galerkin (BUG) integrators \cite{EKDiscretized, GC22Unconventional, CKL22} are stable in the presence of small or zero singular values including for dissipative dynamical systems. These approaches have also been extended to tensor networks \cite{GC2023, GC21Time}. However, these techniques have a first-order temporal accuracy.  Recently, a second-order robust BUG integrator based on the midpoint rule was introduced \cite{CEKL24}.   In \cite{EK19Projection}, a projection method was introduced, in which higher-order Runge-Kutta schemes are applied to the projected equation. In \cite{RDV22}, a similar approach was presented, which is based on the rank-truncation of the time-discrete evolution equations and can achieve high-order temporal accuracy.  

The focus of the above studies has been on developing stable time integration schemes. There are fewer works focused on addressing the issue of computational cost for nonlinear MDEs/TDEs. The issue of computational cost was addressed in  
\cite{MNAdaptive}, in the context of DLRA for MDEs, where a low-rank approximation of $\mathcal{F}(V)$ is constructed using an interpolatory CUR decomposition. The approach presented in \cite{MNAdaptive} requires evaluating $\mathcal{F}(V)$ at $r$ columns and $r$ rows. The columns and rows are selected based on the Discrete Empirical Interpolation Method (DEIM) \cite{SCNonlinear}.  This approach evaluates $\mathcal{F}(V)$ at $r(n_1+n_2)-r^2$ entries, which are the minimum number of entries required to build a rank-$r$ approximation of $\mathcal{F}(V)$ of size $n_1 \times n_2$. A similar technique has also been introduced in \cite{PV23} for  hyper-reduction 
of parametric Hamiltonian dynamical systems. At the time of preparation of this manuscript, a similar interpolatory low-rank approximation \cite{D24} was proposed for tensor train integration which interpolates $\mathcal{F}(V)$ onto the tangent space. All of these techniques are based on DLRA and tangent space projections. We recently extended this approach to DLRA of TDEs in the Tucker tensor form \cite{GB24}. For an overview of the DEIM algorithm for nonlinear reduced order modeling, we refer the reader to  \cite{GB24}.

Recently, a new algorithm for the time integration of nonlinear MDEs on low-rank manifolds was developed that addresses the aforementioned issues \cite{MDOblique}. The algorithm, presented in \cite{MDOblique}, does not utilize tangent space projection. Instead, it applies an interpolatory CUR low-rank approximation to the time-discrete MDE. The algorithm is cost-optimal regardless of the type of nonlinearity of $\mathcal{F}(V)$, evaluating only $r(n_1+n_2)$ entries of the time-discrete MDE at each time instant. Moreover, this algorithm does not rely on the inversion of singular values; it is robust in the presence of small or zero singular values and can achieve high-order temporal accuracy.

In this paper, we introduce a novel algorithm for the time integration of TDEs on low-rank manifolds, addressing the aforementioned challenges. In particular, we build upon our previous work \cite{MDOblique} and develop new time integration schemes for TDEs in tensor train and Tucker tensor forms. The combination of these forms can generate both binary and non-binary tensor tree networks, therefore showcasing the potential applicability of the developed techniques to complex tensor networks. The schematic of the presented algorithm is shown in Figure \ref{fig:overview}. The contributions of this paper are listed below:

(i) We present \texttt{TT-CUR-DEIM (iterative)}, a new tensor train CUR algorithm in which entries are selected based on DEIM \cite{SCNonlinear}. 
The algorithm is similar to the TT-CROSS-maxvol algorithm from \cite{ISTtcross} with the primary difference that the maxvol algorithm is replaced with DEIM.  The DEIM algorithm requires access to exact or approximate singular vectors, which are not available in many applications. Unlike DEIM, the \texttt{TT-CUR-DEIM (iterative)} algorithm does not require access to the singular vectors and can be used as a black box tensor train cross algorithm. If exact or approximate singular vectors are available, no iterations are required. In such cases, we refer to the algorithm simply as \texttt{TT-CUR-DEIM}.

(ii) We propose a novel time integration approach for solving TDEs in low-rank tensor train form, which addresses the challenges highlighted above.  The proposed methodology requires evaluation of $\mathcal{F}(V)$ at $\mathcal{O}(dnr^2)$ entries using the TT-CUR-DEIM algorithm. Similar to the approach presented in \cite{MDOblique},  \texttt{TT-CUR-DEIM}  is designed to be agnostic to the type of nonlinearity encountered (e.g., fractional, exponential, etc.).  Implementing the presented methodology is straightforward, as it requires only the evaluation of the time-discrete FOM at strategically selected entries. We numerically demonstrate that the time integration using the cross approximation is robust in the presence of small singular values and it is extended for high-order explicit time integration schemes, such as the fourth-order Runge-Kutta (RK4). The algorithm is also rank-adaptive, effectively adjusting the rank of the approximation on the fly.

(iii) We extend the above time integration approach to solve TDEs in the low-rank Tucker tensor form, utilizing a recently developed Tucker cross algorithm based on DEIM \cite{GB24}. We follow steps analogous to those for solving TDEs in the tensor train form or solving MDEs in the low-rank form \cite{MDOblique}. This suggests the potential generalizability of the presented framework to more complex tensor networks.

The remainder of the paper is organized as follows. The methodology is discussed in Section 2, demonstrations and results are presented in Section 3, and the conclusions follow in Section 4.

\section{Methodology}
\label{methodology}

\subsection{Preliminaries}
\label{sec:Definition}
We first introduce the notation used in this paper. Vectors are denoted in bold lowercase letters (e.g. $\bm v$),  matrices are denoted by bold uppercase letters (e.g. $\bm V$), and tensors by uppercase letters (e.g. $V$). The entries of tensor $V$ are shown with $V(i_1, \dots, i_d)$.   We denote the mode-$k$ unfolding of tensor $V$ to a matrix with: 
\begin{equation*}
\bm V_{(k)} = [V(i_k,~ i_1 \dots i_{k-1} i_{k+1} \dots i_d)],
\end{equation*}
where $\bm V_{(k)}$ is a matrix of size $n_k \times n_1 \dots n_{k-1} n_{k+1} \dots n_d$.

To extract submatrices, we use integer sets that contain the indices of the submatrix. For example, let $\mathcal I = \big \{i_1^{(\alpha)}\big \}_{\alpha=1}^r$ be a set containing the indices of the rows of $\bm V$, i.e., $\mathcal I$ contains $r$ integer indices such that  $\mathcal I  \subset \{1,2,\dots,n_1 \}$.

If $\mathcal I$ contains \emph{all} indices of rows or columns of $\bm V$, we use a colon.  For example, consider a matrix $\bm V \in \mathbb{R}^{n_1 \times n_2}$ and let $\mathcal I = \big \{i_1^{(\alpha)}\big \}_{\alpha=1}^r$  and $\mathcal J = \big \{i_2^{(\alpha)}\big \}_{\alpha=1}^r$, then $\bm V(\mathcal I,:) \in \mathbb{R}^{r \times n_2}$ is a submatrix that contains $r$ rows of matrix $\bm V$ with the row indices included in the set $\mathcal I$. Similarly, $\bm V(:,\mathcal J) \in \mathbb{R}^{n_1 \times r}$ is a submatrix containing $r$ columns of $\bm V$ and $\bm V(\mathcal I,\mathcal J) \in \mathbb{R}^{r \times r}$ is a submatrix containing the intersection of $r$ rows and columns of $\bm V$. 

Subtensors can be extracted similarly. For example, let $V \in \mathbb{R}^{n_1\times \dots \times n_d}$ and $\mathcal I = \big \{i_1^{(\alpha)}\big \}_{\alpha=1}^r$, then $V(\mathcal I,:) \in \mathbb{R}^{r \times n_2 \times \dots \times n_d}$ is a subtensor of $V$ containing  $r$ indices of the first mode, where we use a single colon to denote all the indices in the remaining $d-1$ modes. We also use tuples to denote the set of indices of multiple modes. For example, let $\mathcal I = \big \{i_1^{(\alpha)}, i_2^{(\alpha)}, i_3^{(\alpha)} \big \}_{\alpha=1}^r$  and $\mathcal J = \big \{i_4^{(\alpha)}, i_5^{(\alpha)}, \dots,i_d^{(\alpha)} \big \}_{\alpha=1}^r$. An example of the set $\mathcal{I}$ for $r=2$ is: $\mathcal I =\{(8,2,9), (5,1,3) \}$. In this case,  $V(\mathcal I,:)$ is a subtensor of $V$ of size $V(\mathcal I,:) \in \mathbb{R}^{r\times n_4 \times n_5 \times \dots \times n_d}$ where the indices of the first three modes are constrained to any of the $r$ members of the set $\mathcal I$  and similarly $V(:,\mathcal J)$ is a subtensor of $V$ of size $V(:,\mathcal J) \in \mathbb{R}^{n_1 \times n_2 \times n_3 \times r}$. 

The symbol $\bigtimes_n$ is used to denote the $n$-mode product. The $n$-mode product of a tensor $V \in \mathbb{R}^{n_1 \times n_2 \times \dots \times n_d}$ with a matrix $\bm B \in \mathbb{R}^{m \times n_n}$ is obtained by $V \bigtimes_n \bm B$ and is of size $n_1 \times \dots \times n_{n-1} \times m \times n_{n+1} \times \dots \times n_d$.

The Frobenius norm of a tensor is shown by $\| V \|_F$ and is defined as:
\begin{equation}\label{eq:FrobNorm}
    \| V\|_F =  \sqrt{\sum_{i_1=1}^{n_1} \sum_{i_2=1}^{n_2} \dots \sum_{i_d=1}^{n_d} V(i_1, \dots i_d)^2}.
\end{equation}

 We use the following notation for computing the rank-$r$ singular value decomposition (SVD) of a matrix $\bm V \in \mathbb{R}^{m\times n}$. For example, $[\bm U, \boldsymbol \Sigma, \bm Y] =\texttt{SVD}(\bm V,r)$ means computing the SVD of $\bm V$ and truncating at rank $r$, where $r<m$ and $r<n$, $\bm U \in \mathbb{R}^{n\times r}$ is the matrix of left singular vectors,  $\boldsymbol \Sigma \in \mathbb{R}^{r \times r}$ is the matrix of singular values and $\bm Y \in \mathbb{R}^{m \times r}$ is the matrix of right singular vectors.  If any of the singular matrices are not needed, the symbol ($\sim$) is used. For example, $[\bm U,\sim, \sim] =\texttt{SVD}(\bm V,r)$ returns only the first $r$ left singular vectors. We use typewriter font to denote algorithms, e.g., \texttt{SVD} and \texttt{DEIM} and \texttt{reshape}. We use \texttt{permute} to rearrange the dimension of a tensor. For example, to rearrrange the dimension of $V(i_1, i_2, \dots, i_d)$ to $V(i_d, \dots, i_2, i_1)$, we use $V \leftarrow \texttt{permute}(V, [d, d-1, \dots, 1])$.

In the following, we define tensor train and Tucker tensor low-rank approximation models.
 \begin{definition}[Tensor train format]
A tensor $\Vhat \in \mathbb{R}^{n_1 \times \dots \times n_d}$ is expressed in tensor train form as
\begin{equation}\label{eq:TTform}
    \Vhat(i_1, \dots, i_d) = \sum_{\alpha_1=1}^{r_1} \dots \sum_{\alpha_{d-1}=1}^{r_{d-1}} G_1(1,i_1,\alpha_1) \cdot G_2(\alpha_1,i_2,\alpha_2)  \cdots G_d(\alpha_{d-1},i_d,1),
\end{equation}
where $G_i \in \mathbb{R}^{r_{i-1}\times n_i \times r_i}$ are the core tensors, which we refer to as tensor train or TT core. Moreover, $r_0=r_d=1$ and $\bm r = (r_1,\dots, r_{d-1})$ are tensor train compression ranks.
 
 \end{definition}
 \begin{definition}[Tucker tensor format]
A tensor $\Vhat \in \mathbb{R}^{n_1 \times \dots \times n_d}$ is expressed in Tucker tensor  form as
\begin{equation}\label{eq:TUform}
    \Vhat(i_1, \dots, i_d) = \sum_{\alpha_1=1}^{r_1} \dots \sum_{\alpha_{d}=1}^{r_{d}} G(\alpha_1,\dots, \alpha_d) \cdot \bm U_1(i_1,\alpha_1)  \cdots \bm U_d(i_d,\alpha_d),
\end{equation}
where $G \in \mathbb{R}^{r_1 \times \dots \times r_d}$ is the core tensor,  $\bm U_i \in \mathbb{R}^{n_i \times r_i}$ are the factor matrices, and $\bm r=(r_1, \dots, r_d)$ are Tucker tensor multi-rank.
 
 \end{definition}
 \begin{definition}[Low-rank tensor manifolds]\label{def:Mr}
 The low-rank tensor manifold $\mathcal{M}_r$ is defined as the set
\begin{equation*}
\mathcal{M}_r = \{\hat{V} \in \mathbb{R}^{n_1 \times \dots \times  n_d}: \ \mbox{rank}(\hat{V}) = \bm r \}, 
\end{equation*}
of tensors of fixed rank $\bm r$. Any member of the set $\mathcal{M}_r$ is denoted by a hat symbol $( \hat{ \ \ } )$, e.g., $\hat{V}$. We use $\mathcal{M}_r$ to generically refer to low-rank manifolds, represented by either Tucker tensors or tensor trains, with the specific manifold being understood from the context. 
\end{definition}

\begin{figure}[t]
\centering
\includegraphics[width=\textwidth]{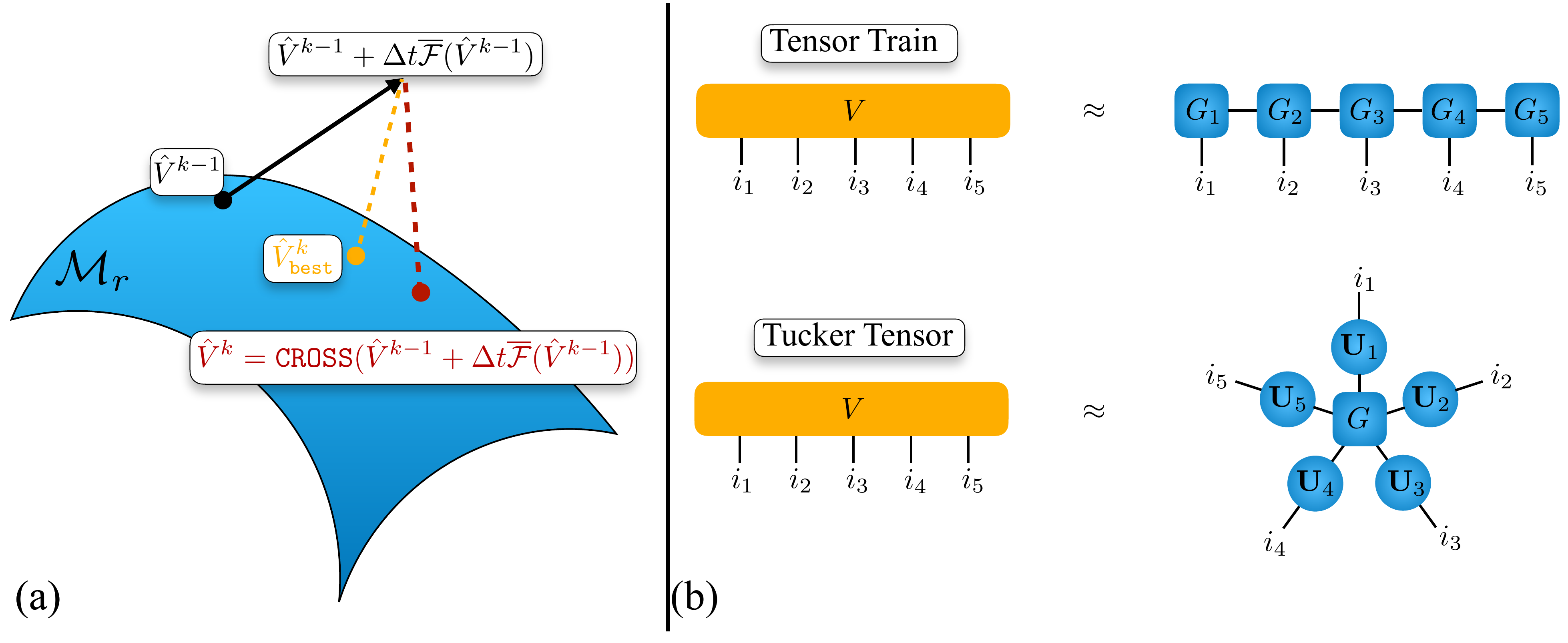}
\caption{(a) Schematic of the cross algorithm for time integration of tensor differential equations on the manifold of low-rank tensors without utilizing tangent space projection. (b) Graphical representation of tensor train and Tucker tensor trees for a five-dimensional tensor $V(i_1,i_2,i_3,i_4,i_5)$.}
\label{fig:overview}
\end{figure}

In the computational complexity analyses, we assume $n \sim \mathcal{O}(n_i)$ and $r \sim \mathcal{O}(r_i)$ for simplicity. Therefore,  the full-dimensional tensor $V \in \mathbb{R}^{n_1 \times \dots \times n_d}$ has $\mathcal{O}(n^d)$ entries, the rank-$r$ Tucker tensor approximation of $V$ reduces the number of entries to $\mathcal{O}(r^d + dnr)$ and   the rank-$r$ tensor train approximation of $V$ reduces the number of entries to $\mathcal{O}(dnr^2)$.

We use the graphical representation of tensors that is common in physics and quantum chemistry \cite{O14}. In this representation, tensors are depicted as nodes of various geometric shapes. See Figure \ref{fig:overview} for a graphical representation of tensor train and Tucker tensor low-rank approximations. The number of legs indicates the order of the tensor. For example, a matrix has two open legs. The interconnecting line between two nodes represents the contraction of the two tensors; it indicates a summation of products over the shared index.

\subsection{Time integration using tangent space projection}
\label{OrthogonalProject}

We consider a general PDE given by:
\begin{equation}\label{eq:ContinousPDE}
    \frac{\partial v (\bm x,t)}{\partial t} = f\bigl( v(\bm x,t) \bigr),
\end{equation}
augmented with appropriate initial and boundary conditions. Here $\bm x = (x_1, x_2, \dots , x_d) \in \mathbb{R}^d$,~ $t$ is time, $f$ is a general nonlinear differential operator, and $d$ is the dimension of the problem. Discretizing the differential operators of  Eq. \ref{eq:ContinousPDE} in $\bm x$ using a method of lines results in the following TDE:
\begin{equation}\label{eq:TDE}
    \frac{d V}{dt} = \mathcal F\bigl( V \bigr),
\end{equation}
where $V(t) \in \mathbb{R}^{n_1 \times n_2 \times \dots \times n_d}$ is the solution tensor and $\mathcal F( \sim )$ is the discrete representation of $f(\ \cdot \ )$. Here, $n_1, n_2, \dots, n_d$ are the number of the discretized degrees of freedom along each mode of the tensor. We refer to Eq. \ref{eq:TDE}  as the full-order model (FOM).

The degrees of freedom of the FOM increase exponentially as $d$ grows. Tensor low-rank approximations can mitigate this issue via low-rank approximation of $V(t)$ by constraining the solution of the above TDE to a manifold of low-rank tensors. In this work, we consider two tensor low-rank approximation schemes: Tucker tensor and tensor train decompositions.  We seek an approximation $\hat{V}(t)$ to the exact solution $V(t)$ where $\hat{V}(t)$ is the constrained solution to low-rank Tucker tensor or tensor train manifold.

Substituting the low-rank approximation $\hat{V}(t)$ into Eq. \ref{eq:TDE} results in a residual equal to:
\begin{equation} \label{residual}
    \mathcal{R}(\dot{\hat{V}}) = \bigl \|  \dot{\hat{V}} - \mathcal F(\hat{V}) \bigr \|_F.
\end{equation}
The low-rank tensor train and Tucker tensor approximation of Eq. \ref{eq:TDE} is obtained by minimization of the above residual with the constraint that $\dot{\hat{V}} \in \mathcal{T}_{\hat{V}} \mathcal{M}_r$, which results in:
\begin{equation} \label{orthoProj}
    \dot{\hat{V}} = \mathcal{P}_{\mathcal{T}_{\hat{V}}}  \mathcal F( \hat{V}), 
\end{equation}
where $\mathcal{T}_{\hat{V}} \mathcal{M}_r$ is the tangent space of the manifold $\mathcal{M}_r$ at $\hat{V} \in \mathcal{M}_r$, and $\mathcal{P}_{\mathcal{T}_{\hat{V}}}$ is the orthogonal projection onto the tangent space $\mathcal{T}_{\hat{V}}$ at $\hat{V}$. The DLRA evolution equations for the Tucker tensor and tensor train models are presented in \cite{KL10} and \cite{CL2013}, respectively. 

The key advantage of the DLRA formalism is that it could be applied to generic TDEs, thereby broadening the application of these techniques beyond the technical developments tailored to solving quantum many-body problems. It also enriches our understanding of solving TDEs on low-rank manifolds by incorporating Riemannian geometry concepts \cite{UV20,SCK24}.   While this formalism has the potential to significantly reduce the computational cost of solving high-dimensional TDEs, there are still several challenges for practical problems of interest - especially for nonlinear TDEs. As outlined in the introduction, computing $F= \mathcal F(V) \in \mathbb{R}^{n_1 \times \dots \times n_d}$ in many nonlinear problems requires $\mathcal{O}(n^d)$ operations, mirroring the computational complexity of solving the FOM. Even in linear TDEs, a substantial number of terms on the right-hand side may result in a large exact rank for $F$. While it is possible to achieve the reductions offered by tensor low-rank approximations in special cases of linear and quadratic nonlinear TDEs, this requires a highly intrusive and meticulous treatment of the right-hand side terms to avoid storing full-rank or high-rank tensors. Moreover, the projection of Eq. \ref{orthoProj} may encounter computational stability issues since the derived DLRA evolution equations from Eq. \ref{orthoProj} include the inverse of auto-correlation matrices. When the tensor solution incorporates small singular values, the auto-correlation matrices become ill-conditioned. This phenomenon is attributed to the curvature of the tensor manifold $\mathcal{M}_r$, which is inversely proportional to the smallest singular value \cite{CLProjector, ARAdaptive}.

\subsection{Time discrete variational principle and cross low-rank approximation}
\label{ObliqueProject}

To address the challenges of the DLRA, we propose a methodology that is inspired by our previous work on matrix low-rank approximation \cite{MDOblique}. In particular, we consider the temporal discretization of the FOM:

\begin{equation} \label{TimeDiscreteTDE}
    V^k = \hat{V}^{k-1} + \Delta t \overline{\mathcal F}( \hat{V}^{k-1}),
\end{equation}
where $\Delta t$ is the time-step size, and $\overline{\mathcal F}$ depends on the temporal discretization that could be a multi-step or a Runge-Kutta-based scheme. For example, the first-order explicit Euler method is simply $\overline{\mathcal F}(\hat{V}^{k-1}) = \mathcal F(\hat{V}^{k-1})$. It should be noted that $\hat{V}^{k-1} + \Delta t \overline{\mathcal F}( \hat{V}^{k-1})$ is calculated based on the solution at the previous time instant being in the low-rank form.  Depending on the type of $\mathcal  F(\ \cdot \ )$, the tensor $V^k$ is generally of higher rank or full rank. For example, if $\mathcal  F(\ \cdot \ )$ has non-polynomial nonlinearity, $\mathcal  F(\Vhat^{k-1})$ is full rank even when $\Vhat^{k-1}$ is low rank. 
In fact,  except for a few rare cases, one forward step of Eq. \ref{TimeDiscreteTDE} moves $V^k$ out of the rank-$r$ manifold. Thus, to keep the solution of the TDE on the low-rank manifold, a rank truncation is required to retract the solution back onto $\mathcal{M}_r$ at every time step. Ideally, we need to approximate $V^k$ with a rank-$r$ tensor $\hat{V}^k$, such that:
\begin{equation} \label{Errortensor}
  \|R^k\|_F = \|V^{k} - \hat{V}^{k}\|_F,
\end{equation}
 is minimized, where $R^k$  is the time-discrete residual due to the low-rank approximation. The best approximation can be found by minimizing the norm of the residual \cite{EK19Projection,CKL22}:
\begin{equation} \label{Residual2}
   \mbox{min} \ \mathcal{J}(\hat{V}^k) = \Bigl \| V^{k} -\hat{V}^{k}  \Bigr \|_F^2, \quad \mbox{such that} \quad \hat{V}^k \in \mathcal{M}_r.
\end{equation}
The best solution is shown in Figure \ref{fig:overview} as $\hat{V}^k_{\texttt{best}}$. 
We refer to the minimization described above as the time-discrete variational principle.

Solving this minimization problem is quite costly, not only because it requires an iterative algorithm but also because it involves forming the full or high-rank tensor $V^k$. Avoiding the computational cost of calculating $V^k$—both in terms of memory requirements and FLOP — is precisely why we pursue low-rank approximation in the first place.

Another approach to finding a near-optimal $\hat{V}^k$ involves using tensor train SVD (TT-SVD) for tensor train decomposition, or higher-order SVD (HOSVD) for Tucker tensor decomposition. Although neither of these approaches is iterative, any SVD-based methods, including randomized algorithms, require computing all entries of $V^k$ when $V^k$ is of full rank. This is particularly necessary, for instance, for non-polynomial nonlinear $\mathcal{F}(\ \cdot \ )$.

Our approach is to construct the low-rank approximation tensor $\hat{V}^{k}$ using a cross approximation as shown below:
\begin{equation} \label{TDE-CROSS}
    \hat{V}^{k} = \texttt{CROSS}(\hat{V}^{k-1} + \Delta t \overline{\mathcal  F}(\hat{V}^{k-1})).
\end{equation}
The key advantage of the above approach is that $\hat{V}^{k}$ can be constructed by evaluating $V^k = \hat{V}^{k-1} + \Delta t \overline{\mathcal{F}}(\hat{V}^{k-1})$ at judiciously chosen entries, thereby avoiding the formation of full-rank tensors. However, the accuracy of cross algorithms depends on which entries are chosen. In the following section, we present a DEIM cross algorithm for constructing tensor train low-rank approximation that in practice is nearly as accurate as TT-SVD. In the case of the Tucker tensor model, we will utilize the recently developed cross-algorithm \cite{GB24}.

\subsection{Cross low-rank approximations based on DEIM}
\label{TT-CUR-DEIM}
In this section, we first review the CUR-DEIM algorithm for matrix low-rank approximation. We then present a new cross algorithm for tensor train low-rank approximation based on CUR-DEIM.

\subsubsection{CUR-DEIM algorithm for matrix low-rank approximation}
We first present the CUR algorithm based on DEIM. This algorithm was utilized in \cite{MDOblique, MNAdaptive} to solve nonlinear matrix differential equations on low-rank manifolds. Let $\bm V \in \mathbb{R}^{n_1 \times n_2}$ be a matrix and let $\mathcal I = \{ i_1^{(\alpha)}\}$ and $\mathcal J = \{ i_2^{(\alpha)}\}$ for $\alpha = 1, \dots, r$ be a set of row and column indices, respectively. Therefore, $\bm R = \bm V(\mathcal I,:) \in \mathbb{R}^{r \times n_2} $ and $\bm C = \bm V(:,\mathcal J) \in \mathbb{R}^{n_1 \times r} $ are the row and column submatrices of $\bm V$.    A CUR rank-$r$ approximation of $\bm V$ can be obtained by $\hat{\bm V} = \bm C \tilde{\bm U} \bm R$,
where $\tilde{\bm U} \in \mathbb R^{r \times r}$. The accuracy of the above low-rank approximation depends on $\mathcal I$ and $\mathcal J$ and on how $\tilde{\bm U}$ is computed. The best matrix $\tilde{\bm U}$ can be computed via the orthogonal projection of $\bm V$ onto the space spanned by $\bm C$ and $\bm R$, i.e., $\tilde{\bm U} = (\bm C^T \bm C)^{-1} \bm C^T \bm V \bm R^T (\bm R \bm R^T)^{-1}$ \cite{SE16}. However, orthogonal projection requires all elements of the matrix $\bm V$, which is undesirable if the elements of $\bm V$ are costly to compute. A remedy to this problem is the use of interpolatory projection, where $\tilde{\bm U} = \bm V^{-1}(\mathcal I, \mathcal J)$. The resulting CUR low-rank approximation is given by:
\begin{equation}\label{eq:CUR_unst}
\hat{\bm V} = \bm V(:,\mathcal J) \bm V^{-1}(\mathcal I, \mathcal J) \bm V(\mathcal I,:).
\end{equation}
The CUR algorithm, based on interpolatory projectors, utilizes only $s=r(n_1 + n_2) - r^2$ elements of $\bm V$ to build $\hat{\bm V}$. This is the minimum number of elements required to build a rank-$r$ approximation of $\bm V$ and for $n=n_1=n_2$ and $r \ll n$, $s\approx 2rn$, which is by a factor of $n^2/2rn=n/2r$ smaller than the size of $\bm V$.  See \cite{MDOblique} for the definition of interpolatory projectors.  It is straightforward to show that CUR-DEIM is a matrix interpolation algorithm because $\bm V(\mathcal I,:) = \hat{\bm V}(\mathcal I,:)$ and $\bm V(:,\mathcal J) = \hat{\bm V}(:,\mathcal J)$.

As mentioned above, the choice of $\mathcal I$ and $\mathcal J$ affects the low-rank approximation error: $\bm E = \bm V - \hat{\bm V}$. To this end, we use the DEIM algorithm to compute $\mathcal I$ and  $\mathcal J$. The DEIM algorithm requires access to the column and row singular vectors of $\bm V$ \cite{CS10} and it is presented in Algorithm \ref{alg:DEIM} in  Appendix A for convenience.   Let $\bm U \in \mathbb{R}^{n_1 \times r}$ and $\bm Y \in \mathbb{R}^{n_2 \times r}$ be exact or approximate dominant left and right singular vectors of $\bm V$, respectively. Then $\mathcal I = \texttt{DEIM}(\bm U)$ and $\mathcal J = \texttt{DEIM}(\bm Y)$. We refer to the CUR algorithm based on DEIM as CUR-DEIM. According to \cite[Theorem 2.8]{MDOblique}, the CUR-DEIM low-rank approximation error is bounded by:
\begin{equation}
    \|\bm V -\hat{\bm V} \|_2 \leq c \hat{\sigma}_{r+1},
\end{equation}
where $\hat{\sigma}_{r+1}$ is the orthogonal projection error:  $\hat{\sigma}_{r+1} = \mbox{max}\{ \|(\bm I - \bm U \bm U^T) \bm V  \|_2, \| \bm V(\bm I - \bm Y \bm Y^T)   \|_2 \}$, and $\| \cdot \|_2$ is the second matrix norm. When $\bm U$ and $\bm Y$ are the exact singular vectors, $\hat{\sigma}_{r+1}$ is the $(r+1)$-th singular value of $\bm V$, i.e., the lowest error for a rank-$r$ approximation. The amplification factor $(c\geq 1)$ depends on the condition number of two matrices: $c = \mbox{min} \{ \eta_r (1+\eta_c), \eta_c(1+\eta_r)\}$, where $\eta_r =\|\bm U^{-1}(\mathcal I,:) \|_2 $ and $\eta_c =\|\bm Y^{-1}(\mathcal J,:) \|_2 $. An important result of \cite[Lemma 3.2]{CS10} is that $\eta_r$ (or $\eta_c$) is bounded by:
\begin{equation}
\eta_r \leq \|(1+ \sqrt{2n_1})^{r-1} \|\bm u_1 \|^{-1}_{\infty},
\end{equation}
where $\bm u_1$ is the first left singular vector, i.e., the first column of $\bm U$. A similar bound can be shown for $\eta_c$. The above result shows that $\eta_r$ (or $\eta_c$) remains bounded irrespective of the singular values- thereby guaranteeing the DEIM remains well-conditioned as $r$ increases. The above error bound is quite pessimistic, and in practice, $\eta_r$ and $\eta_c$, which can be computed cheaply, are small. In fact,  DEIM computes the interpolation points based on a greedy algorithm such that $\eta_r$ and $\eta_c$ are minimized.

 The DEIM algorithm requires either exact or approximate left and right singular vectors of the matrix $\bm{V}$ to determine $\mathcal{I}$ and $\mathcal{J}$. However, computing the singular vectors of $\bm{V}$ necessitates access to all its elements, which is undesirable when such access is costly. This is often the case, for example, when the elements of $\bm{V}$ must be computed, as with MDEs. This issue can be resolved by using CUR algorithms for solving MDE/TDE, where $\mathcal{I}$ and $\mathcal{J}$ are computed based on the singular vectors from the previous time step. However, this issue must be addressed if we are to use CUR-DEIM as a black-box low-rank approximation algorithm.

We present a modified CUR algorithm that addresses the aforementioned practical issue. It is possible to relax the requirement from needing access to both left and right singular vectors to just one set of vectors, either left or right. For example, assume the column indices ($\mathcal J$) are obtained by applying DEIM on the right singular vectors $\bm Y$. Applying SVD on the selected columns results in:
\begin{equation}\label{eq:SVD_DEIM}
    [\bm U, \boldsymbol \Sigma, \sim] = \texttt{SVD}(\bm V(:,\mathcal J),r).
\end{equation}
 The DEIM algorithm can now be applied to $\bm U$ to compute the row indices:  $\mathcal I = \texttt{DEIM}(\bm U)$.    Interpolating the columns of $\bm V$ onto $\bm U$ using the selected rows results in: 
\begin{equation}\label{eq:CUR_st}
\hat{\bm V} = \bm G \bm R,
\end{equation}
where $\bm G = \bm U \bm U^{-1}(\mathcal I,:)$ and $\bm R= \bm V(\mathcal I,:)$. Therefore, $\mathcal{I}$ does not need to be provided to the CUR algorithm; it can be computed internally. However,  $\mathcal J$ still needs to be provided.  To address this issue, we present an iterative CUR-DEIM algorithm later in this section that does not require access to the singular vectors of $\bm{V}$. We introduce this iterative algorithm for tensor train low-rank approximation, of which matrix low-rank approximation is a special case, i.e., when $d=2$. Another advantage of computing CUR using Eq. \ref{eq:CUR_st} is that it is numerically stable. A similar solution based on QR decomposition and the maxvol algorithm was presented in \cite{TT11}, which also results in a stable CUR algorithm.   We also present a rank-adaptive CUR algorithm in which the rank is determined adaptively based on an accuracy threshold.

\subsubsection{CUR-DEIM algorithm for tensor train low-rank approximation}\label{sec:TT-CUR-DIEM}
We present a CUR-DEIM algorithm for constructing a low-rank tensor train approximation. We refer to this algorithm as $\texttt{TT-CUR-DEIM}$. Let $V \in \mathbb{R}^{n_1 \times \dots \times  n_d}$ be the target tensor and let $\bm V_{(1)} =[V(i_1; i_2, \dots, i_d)] \in \mathbb{R}^{n_1 \times n_2 n_3 \dots n_d }$ be the mode-1 unfolding matrix of $V$. Therefore, each column of $\bm V_{(1)}$ is a vector of size $n_1 \times 1$ and each row of  $\bm V_{(1)}$ can be reshaped to a $(d-1)$-dimensional tensor of size $n_2\times \dots \times n_d$.

To produce the first tensor train core, we apply \texttt{CUR-DEIM} on matrix $\bm V_{(1)}$. The index of each column of $\bm V_{(1)}$ can be mapped uniquely to a collection of $(d-1)$-tuple index of $\big \{i_2,\dots, i_d \big \}$.  Let $\mathcal J_1 = \big\{i^{(\alpha_1)}_2, \dots, i^{(\alpha_1)}_d\big\}_{\alpha_1=1}^{r_1}$ be a set of $(d-1)$-tuples and $\alpha_1=1, \dots, r_1$. Assuming $\mathcal J_1$ is known, $\bm U_1 \in \mathbb{R}^{n_1 \times r_1}$ can be  computed following the CUR algorithm as follows:
\begin{equation*}
    [\bm U_1, \sim, \sim ] = \texttt{SVD}(V(: ,\mathcal J_1),r_1),
\end{equation*}
where $V(: ,\mathcal J_1) \in \mathbb{R}^{n_1\times r_1}$ is the submatrix that contains $r_1$ columns of $\bm V_{(1)}$. Next, $\mathcal I_1$ is calculated:
\begin{equation*}
    \mathcal I_1 = \texttt{DEIM}(\bm U_1), \quad \mbox{where} \quad \mathcal I_1 = \big\{i_1^{(\alpha_1)}\big\}, \quad \alpha_1=1, \dots, r_1.
\end{equation*}
Subsequently, the matrix $\bm G_1 = \bm U_1 \bm U_1^{-1}(\mathcal I_1,:)$ is calculated.  The matrix $\bm G_1 \in \mathbb{R}^{ n_1 \times r_1}$ can be reshaped to a third-order tensor of size $G_1 \in \mathbb{R}^{1 \times n_1 \times r_1}$. The tensor $G_1$ is the first TT core of $V$.

Applying CUR on matrix $\bm V_{(1)}$ results in a rank-$r_1$ approximation as shown below:
\begin{equation*}
    \hat{\bm V}_{(1)} = \bm G_1 \bm R_1,
\end{equation*}
where $\bm R_1 = \bm V_{(1)} (\mathcal I_1,:) \in \mathbb{R}^{r_1 \times n_2 \dots n_d}$. Note that $\bm{R}_1$ is generally a very large matrix; as we demonstrate below, it is never explicitly formed or stored. 

\begin{figure}[t]
\centering
\includegraphics[width=\textwidth]{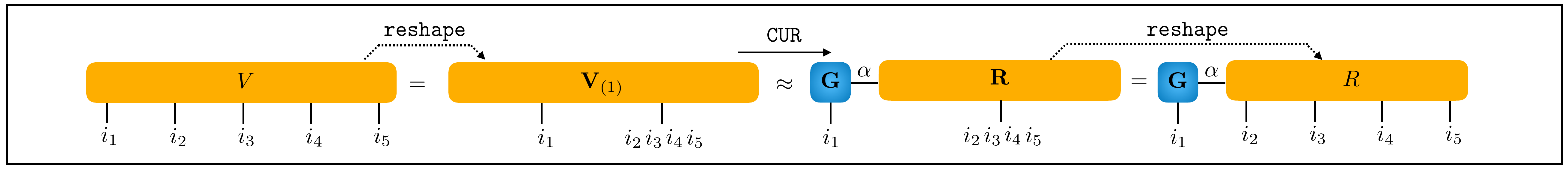}
\caption{Graphical representation of \texttt{CUR-DEIM} algorithm applied to the  first unfolding of a fifth-order tensor $V(i_1,i_2,i_3,i_4,i_5)$, where $V(i_1,i_2,i_3,i_4,i_5) \approx \sum_{\alpha=1}^r \mathbf G(i_1,\alpha) R(\alpha,i_2,i_3,i_4,i_5) $. In this representation,  the single index belongs to the integer set of $[1,2,\dots, n_k]$ and the combined indices $i_2i_3 i_4 i_5$ belongs to the integer set of $[1,2,\dots, n_2 n_3 n_4 n_5]$. }
\label{fig:CUR-DEIM}
\end{figure}

The above steps  can be expressed  as:
\begin{equation}
[\bm G_1, \mathcal I_1] = \texttt{CUR-DEIM}(V,\mathcal J_1),
\end{equation}
where  \texttt{CUR-DEIM}  is presented in Algorithm \ref{alg:CUR-DEIM}. In Figure \ref{fig:CUR-DEIM}, the graphical representation of the \texttt{CUR-DEIM} algorithm is shown.  Note that the \texttt{CUR-DEIM} algorithm evaluates $V$ at a small number of elements, on the order of $\mathcal{O}(n)$, and $V$, as an input to this algorithm, should be interpreted as a function that returns the values of $V$ at the requested entries.

\begin{algorithm}[t]
\fontsize{9pt}{9.5pt}\selectfont
\SetAlgoLined
\KwIn{$V$: function handle that returns $V(i_1, i_2, \dots, i_m), \quad$  $\mathcal{J}= \{i_2^{(\alpha)}, \dots, i_m^{(\alpha)} \}, \quad \alpha=1, \dots, r$}
\KwOut{$ \bm G, ~ \mathcal I$}

$[\bm U, \sim, \sim ] = \texttt{SVD}(V(:, \mathcal{J}),r)$ \hspace{1.95cm} $\rhd$ compute the SVD of the selected columns\;

$\mathcal I =\texttt{DEIM}(\bm U)$ \hspace{3.84cm} $\rhd$ compute the DEIM indices\;

$\bm G =\bm U \bm U^{-1}(\mathcal I,:)$ \hspace{3.4cm} $\rhd$ compute the column factorized matrix\;

\caption{\texttt{CUR-DEIM} algorithm for unfolded tensors.}
\label{alg:CUR-DEIM}
\end{algorithm}

In the next step, we show how the \texttt{CUR-DEIM} algorithm can be applied to $\bm R_1$ to extract the second TT core without explicitly forming $\bm R_1$. The matrix $\bm R_1$ can be reshaped into a $(d-1)$-dimensional tensor of size $R_1 \in \mathbb{R}^{r_1 n_2 \times \dots \times n_d}$ where the indices $\alpha_1$ and $i_2$ are combined into a long index. It is easy to verify that $R_1$ is a subtensor of $V$ as shown below:
\begin{equation}\label{eq:aux_1}
R_1(\alpha_1 i_2, \dots, i_d ) = V(i_1^{(\alpha_1)},i_2, \dots, i_d ), \quad \mbox{where} \quad i_1^{(\alpha_1)} \in \mathcal I_1, \quad \alpha_1=1, \dots, r_1.
\end{equation}

To obtain the second core, we follow the steps analogous to those used for the first core, where $V$ is replaced with  $R_1$. To this end, let $\mathcal J_2 = \big \{i^{(\alpha_2)}_3, \dots, i^{(\alpha_2)}_d\big \}_{\alpha_2=1}^{r_2}$ be a set of $(d-2)$-tuples. Similar to the first core, we assume that $\mathcal J_2$ is known.  
Applying the \texttt{CUR-DEIM} algorithm to $R_1$ results in:
\begin{equation}
[\bm G_2, \mathcal I_2] = \texttt{CUR-DEIM}(R_1,\mathcal J_2),
\end{equation}
where $\bm G_2 \in \mathbb{R}^{r_1 n_2 \times r_2}$. The matrix $\bm G_2$ is reshaped to the tensor $G_2 \in \mathbb{R}^{r_1 \times n_2 \times r_2}$. The tensor $G_2$ is the second TT core. The set $\mathcal I_2$ contains the  $r_2$ DEIM-selected indices that are a subset of the row indices of $\bm G_2$, which is the integer set of $\{1,2 \dots, r_1 n_2\}$. Given the relationship between $V$ and $R_1$ by Eq. \ref{eq:aux_1}, any row index of $\bm G_2$ can be uniquely represented with the index pair of $\{i_1^{(\alpha_1)},i_2 \}$, where $i_1^{(\alpha_1)} \in \mathcal I_1$ and $i_2 \in \{1,2, \dots, n_2\}$. Therefore, the index set $\mathcal I_2$ can be represented as:
\begin{equation}
\mathcal I_2 = \big\{i_1^{(\alpha_2)},i_2^{(\alpha_2)}\big \}_{\alpha_2=1}^{r_2},
\end{equation}
where $\big \{i_1^{(\alpha_2)} \big \} \subset \mathcal I_1 $.

To obtain the third TT core, the same steps are applied to a $(d-2)-$dimensional tensor $R_2 \in \mathbb{R}^{r_2 n_3 \times n_4 \times \dots \times n_d}$, where
\begin{equation}
    R_2(\alpha_2 i_3,i_4, \dots, i_d) = V(i_1^{(\alpha_2)},i_2^{(\alpha_2)},i_3, \dots, i_d).
\end{equation}
Assuming $\mathcal J_3 = \big \{i^{(\alpha_3)}_4, \dots, i^{(\alpha_3)}_d\big \}_{\alpha_3=1}^{r_3}$ is known, the \texttt{CUR-DEIM} algorithm is applied to $R_2$ to compute $G_3$ and $\mathcal I_3$. These steps are analogously applied to generate the subsequent TT cores. The last step of this algorithm involves applying the CUR-DEIM algorithm to 
\begin{equation}
    R_{d-2}(\alpha_{d-2} i_{d-1}, i_d) = V(i_1^{(\alpha_{d-2})},\dots, i_{d-2}^{(\alpha_{d-2})}, i_{d-1}, i_d).
\end{equation}
Assuming $\mathcal J_{d-1} =\big \{i_{d-1}^{\alpha_{d-1}} \big \}_{\alpha_{d-1}=1}^{r_{d-1}}$ in known $G_{d-1}$ and $\mathcal I_{d-1}$ are computed as follows:
\begin{equation}
[\bm G_{d-1}, \mathcal I_{d-1}] = \texttt{CUR-DEIM}(R_{d-2},\mathcal J_{d-1}),
\end{equation}
where $G_{d-1} \in \mathbb{R}^{r_{d-2}\times n_{d-1} \times r_{d-1}}$. The last TT core ($G_d$) can be computed by explicitly forming the rows of the CUR algorithm applied to $R_{d-2}$, which is a matrix of size $r_{d-2}n_{d-1}\times n_d$.  Therefore, $\bm G_d = R_{d-2}(\mathcal I_{d-1},:) \in \mathbb{R}^{r_{d-1}\times n_d}$, where $\mathcal I_{d-1} \subset \{1,2, \dots, r_{d-2}n_{d-1}\} $ is used.  The matrix $\bm G_{d-1}$ is then reshaped to the  TT core tensor of size $G_d \in \mathbb{R}^{r_{d-1}\times n_d \times 1}$. 

\begin{figure}[t]
\centering
\includegraphics[width=\textwidth]{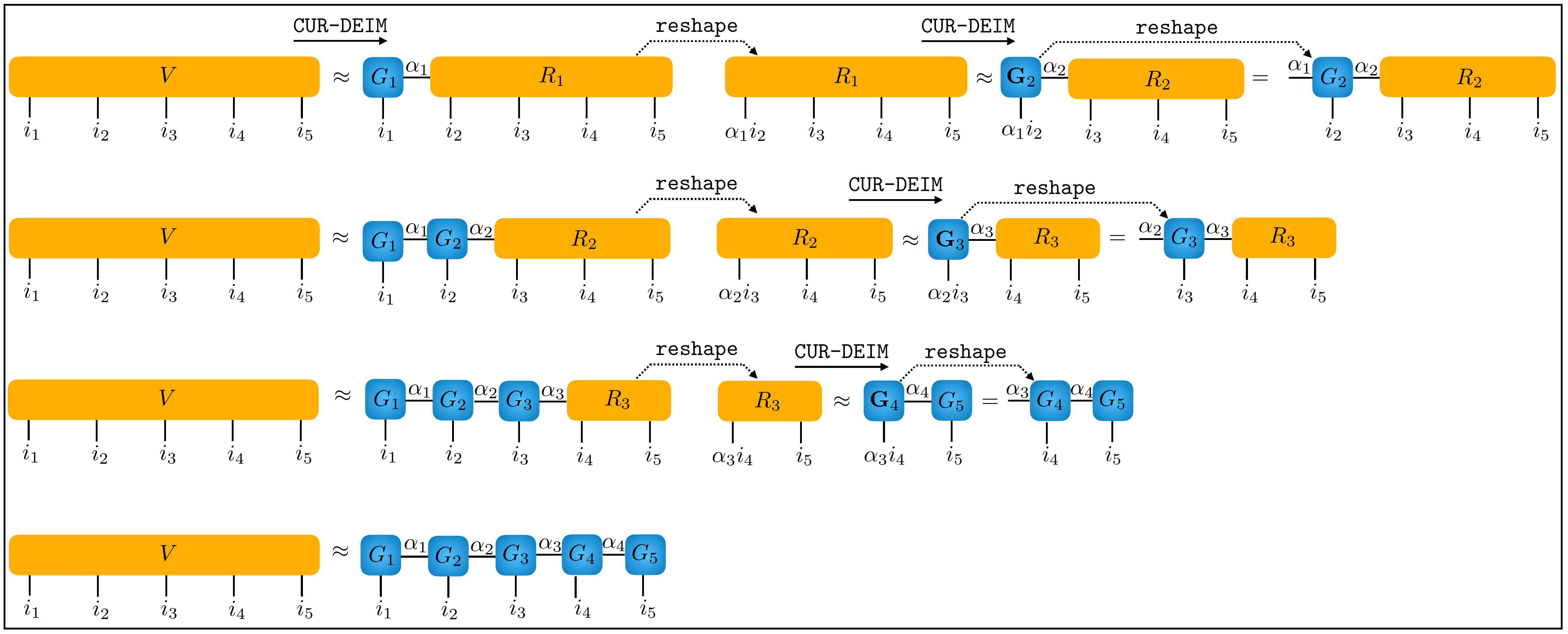}
\caption{Graphical representation of \texttt{TT-CUR-DEIM} algorithm for a fifth-order tensor $V(i_1,i_2,i_3,i_4,i_5)$.}
\label{fig:TTscheme}
\end{figure}

\begin{figure}[t]
\centering
\includegraphics[width=.8\textwidth]{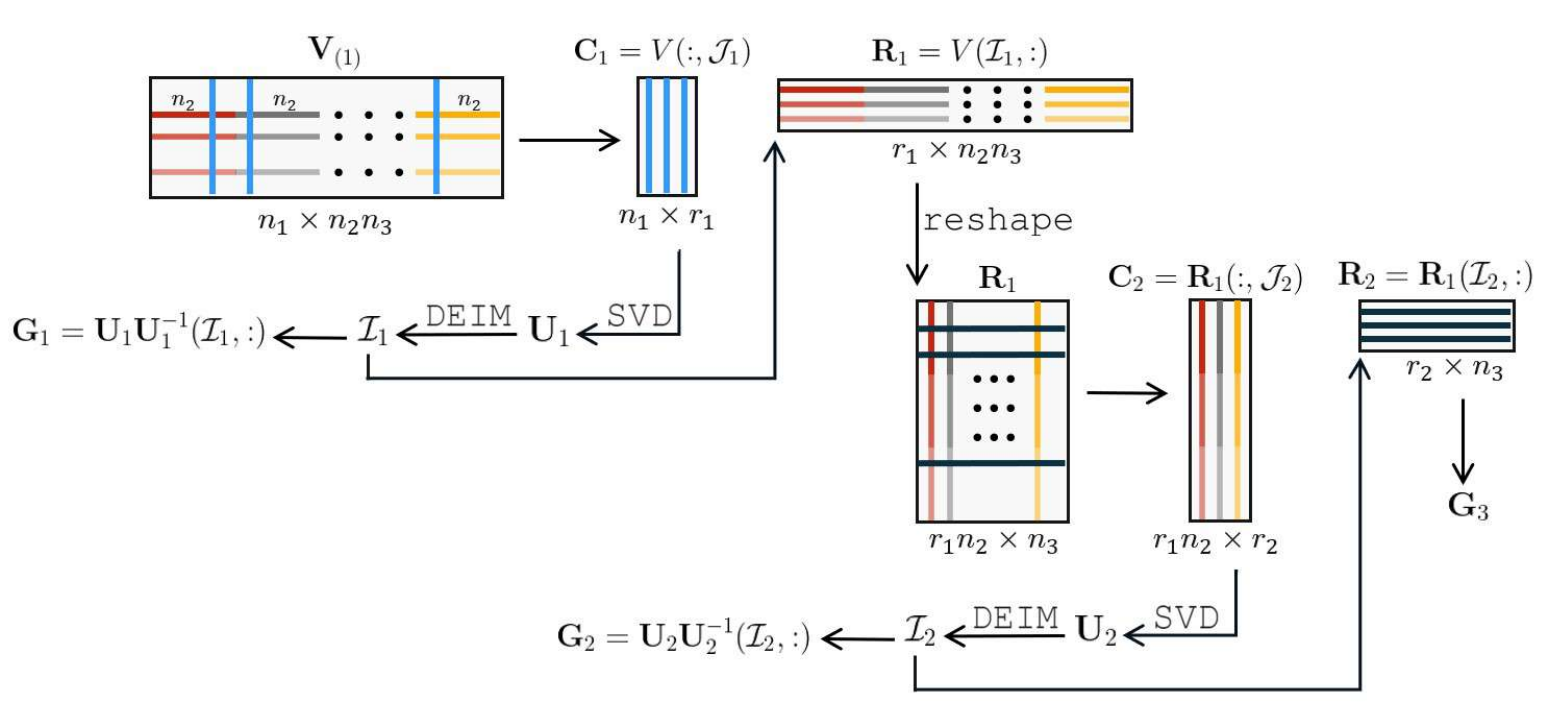}
\caption{Schematic of the \texttt{TT-CUR-DEIM} algorithm for a three-dimensional tensor $V$ of size $n_1 \times n_2 \times n_3$. Note that $\bm R_1$ is not computed or stored and $\bm C_2$ and $\bm R_2$ are obtained directly from $V_{(1)}$.}
\label{fig:TT-CUR-DEIM Scheme}
\end{figure}

The above steps are presented in Algorithm \ref{alg:TT_CUR_DEIM} and can be expressed as:
\begin{equation}
    [G_1,\dots, G_d; \mathcal I_1, \dots, \mathcal I_d] = \texttt{TT-CUR-DEIM}(V,\mathcal J_1,\mathcal J_2, \dots, \mathcal J_{d-1}),
\end{equation}
where $\hat{V}$ is obtained in the tensor train format as shown in Eq. \ref{eq:TTform}.
The graphical representation of the \texttt{TT-CUR-DEIM} is shown for a fifth-order tensor in Figure \ref{fig:TTscheme}. A schematic of the selected entries in the \texttt{TT-CUR-DEIM} algorithm is shown in Figure \ref{fig:TT-CUR-DEIM Scheme} for a third-order tensor. 

We refer to $\mathcal{I}_1, \mathcal{I}_2, \dots, \mathcal{I}_{d-1}$ as left indices and $\mathcal{J}_1, \mathcal{J}_2, \dots, \mathcal{J}_{d-1}$ as right indices. The above algorithm amounts to a \emph{left-to-right} sweep. It is possible to perform a \emph{right-to-left} sweep, where the left indices are the input to the algorithm and the right indices are computed as output. 

In the preceding steps, we assumed the availability of a good choice of indices $\mathcal{J}_1, \mathcal{J}_2, \dots, \mathcal{J}_{d-1}$. A near-optimal choice of these indices is accessible when solving TDEs on low-rank manifolds, where the indices used are those computed from the previous time step. However, in general, such a good choice of indices is not readily available. A relevant example in solving TDEs is obtaining the initial condition in a low-rank form. For many PDEs, the initial condition is full-rank, necessitating its low-rank approximation. Once again, the use of TT-SVD is impractical for large tensors due to memory and FLOP costs and another remedy must be sought.

In the following, we present \texttt{TT-CUR-DEIM (iterative)}, where the indices are determined iteratively. These steps are explained below: (i) the right indices $\mathcal{J}_1, \mathcal{J}_2, \dots, \mathcal{J}_{d-1}$ are randomly initialized; (ii) given the right indices, the \texttt{TT-CUR-DEIM} algorithm is executed in a left-to-right sweep, generating the left indices $\mathcal{I}_1, \mathcal{I}_2, \dots, \mathcal{I}_{d-1}$; (iii) the \texttt{TT-CUR-DEIM} algorithm is executed in a right-to-left sweep, taking the left indices as input and computing the right indices as output; (iv) steps (ii) and (iii) are executed iteratively until convergence is achieved. The maximum number of iterations can be determined by the convergence of the singular values of the core tensors.   In practice, a small number of iterations (less than 10) are needed. These steps are presented in Algorithm \ref{alg:TT_CUR_DEIM_it}.

\let\oldnl\nl
\newcommand{\nonl}{\renewcommand{\nl}{\let\nl\oldnl}}

\begin{algorithm}[t]
\begingroup
\fontsize{9pt}{9pt}\selectfont
\caption{\small{\texttt{TT-CUR-DEIM} algorithm for tensor train low-rank approximation.}}
\label{alg:TT_CUR_DEIM}
\SetAlgoLined

\KwIn{$V$:
function handle that returns $V(i_1, i_2, \dots, i_d)$ \\

\hspace{10mm} $\mathcal{J}_1, \mathcal{J}_2, \dots,  \mathcal{J}_{d-1}$: right indices
} \vspace{0.5mm}

\KwOut{$G_1, G_2, \dots, G_d$: core tensors of the tensor train decomposition\\

\hspace{14mm} $\mathcal{I}_1, \mathcal{I}_2, \dots,  \mathcal{I}_{d-1}$: left indices } \vspace{2mm}

$[\bm G_1, \mathcal{I}_1] = \texttt{CUR-DEIM}(V, \mathcal{J}_1)$ \vspace{0.5 mm}

$R_1 = V(\mathcal{I}_1, :)$ \vspace{1mm}

\For{z = 2~:~d-1}{

$[\bm G_z, \mathcal{I}_z] = \texttt{CUR-DEIM} (R_{z-1}, \mathcal{J}_z)$ \vspace{0.5mm}

$R_z = V(\mathcal{I}_z, :)$ \vspace{0.5mm}

}

$ \bm G_d = R_{d-1}$ \vspace{0.5mm}

\For{z = 1~:~d}{
$G_z =\texttt{reshape} (\bm G_z,[r_{z-1},n_z,r_z])$
}
\endgroup
\end{algorithm}

\begin{algorithm}[t]
\begingroup
\fontsize{9pt}{9pt}\selectfont
\caption{\small{\texttt{TT-CUR-DEIM (iterative)} algorithm for tensor train low-rank approximation.}}
\label{alg:TT_CUR_DEIM_it}
\SetAlgoLined

\KwIn{$V$:
function handle that returns $V(i_1, i_2, \dots, i_d)$ \\

\hspace{10mm} $\mathcal{J}_1, \mathcal{J}_2, \dots,  \mathcal{J}_{d-1}$: randomly initialized right indices.
} \vspace{0.5mm}

\KwOut{$G_1, G_2, \dots, G_d$: core tensors of the tensor train decomposition \\

\hspace{14mm} $\mathcal{I}_1, \mathcal{I}_2, \dots,  \mathcal{I}_{d-1}$: left indices }\vspace{2mm}

\For{$z = 1: ~$max iteration}{
[$G_1, \dots, G_d; ~ \mathcal{I}_1, \dots, \mathcal{I}_{d-1}$] = \texttt{TT-CUR-DEIM}($V,\mathcal J_1, \dots, \mathcal J_{d-1})$ \vspace{0.5 mm}

$\mathcal{J}_1 \leftarrow \mathcal{I}_{d-1};~~ \mathcal{J}_2 \leftarrow \mathcal{I}_{d-2}; ~\dots ~; ~ \mathcal{J}_{d-1} \leftarrow \mathcal{I}_1$ \vspace{0.5 mm}

$V \leftarrow \texttt{permute}\bigl(V, [d, d-1, \dots , 1] \bigr)$ \vspace{0.5 mm}

}
\endgroup
\end{algorithm}

\subsection{Runge-Kutta temporal integration}
\label{high-orderTime}
In this section, we explain how the developed cross algorithms can be used for the efficient time integration of TDEs. We focus on three subjects: (i) evaluation of the time-discrete TDE at DEIM-selected entries; (ii) Runge-Kutta time integration; and (iii) how the sampling indices $\mathcal J_1, \mathcal J_2, \dots, \mathcal J_{d-1}$ are determined without performing any iterations. We address these subjects for tensor train low-rank approximations as the analogous procedure can be applied to the Tucker tensor model.  All of these issues were also faced and addressed for solving MDEs \cite{MDOblique}.  

We first consider the explicit Euler time scheme, where $\overline{\mathcal F}(\hat{V}^{k-1}) 
 \equiv \mathcal F(\hat{V}^{k-1})$. We explain how the time-discrete equation $V^k = \hat{V}^{k-1} + \Delta t \mathcal F(\hat{V}^{k-1})$ can be evaluated at the DEIM-selected entries.  Since $\hat{V}^{k-1}$ is in the low-rank form, its evaluation at any set of sparse entries can be done in a straightforward manner and efficiently. We turn our focus to the sparse evaluation of $F = \mathcal F(\hat{V}^{k-1})$, where $F \in \mathbb{R}^{n_1 \times n_2 \times \dots \times  n_d}$ is a tensor of the same size as $\hat{V}$ and it is obtained by applying $\mathcal F ( \cdot)$ to $\hat{V}^{k-1}$. For the sake of simplicity in the explanation, we assume that $\mathcal F( \cdot )$ is obtained by the discretization of  $f(v) = \partial^2 v/\partial x_1^2 + \partial^2 v/\partial x_2^2 + \dots +\partial^2 v/\partial x_d^2$. Let us consider a case where $\mathcal F( \cdot)$ is obtained by the second-order central finite difference discretization of $f(\cdot)$ on a uniform grid with the grid spacing of $\Delta x_k$ in the $x_k$ direction. In the following, we explain how to compute the first TT core. Computing the other TT cores can be done similarly.

 Obtaining the first core of $\hat{V}^k$ requires evaluating $F(:,\mathcal J_1)$.  However, evaluating $F(:,\mathcal J_1)$  requires access to the values of $\hat{V}^{k-1}(:,\mathcal J_1)$ as well as other entries. To make this clear, consider the discretization of the second term ($ \partial^2 v/\partial x_2^2$) at $i_1 \in [1,2 \dots, n_1]$ and $\mathcal J_1 = \big \{i^{(\alpha_1)}_2,\dots, i^{(\alpha_1)}_d \big \}_{\alpha_1=1}^{r_1}$ entries:
\begin{equation*}
\frac{\hat{V}^{k-1}(i_1,i^{(\alpha_1)}_2+1, \dots, i^{(\alpha_1)}_d)- 2\hat{V}^{k-1}(i_1,i^{(\alpha_1)}_2, \dots, i^{(\alpha_1)}_d)+\hat{V}^{k-1}(i_1,i^{(\alpha_1)}_2-1, \dots, i^{(\alpha_1)}_d)}{\Delta x_2^2}.
\end{equation*}
Computing the above expression requires the values of  $\hat{V}^{k-1}$ at $\big \{i^{(\alpha_1)}_2-1, i^{(\alpha_1)}_2+1 \big \}_{\alpha_1=1}^{r_1}$ entries. Therefore, for the above example,  computing $F(:,\mathcal J_1)$ requires the values of  $\hat{V}^{k-1}$ at an auxiliary set denoted with: $\mathcal J_{1_{(a)}} = \bigcup_{k=2}^d  \big \{i^{(\alpha_1)}_k-1, i^{(\alpha_1)}_k+1 \big \}_{\alpha_1=1}^{r_1}$, where $\bigcup$ denotes the union of these indices. Therefore, 
\begin{equation*}
   V^k(:,\mathcal J_1) = \hat{V}^{k-1}(:,\mathcal J_1) + \Delta t \mathcal F(\hat{V}^{k-1}(:,\overline{\mathcal J}_1)),
\end{equation*}
where $\overline{\mathcal J}_1$ is the union of DEIM-selected entries and the dependent indices:  $\overline{\mathcal J}_1 = \mathcal J_1 \bigcup \mathcal J_{1_{(a)}}$. Since  $\hat{V}^{k-1}$ is in the tensor train form, $\hat{V}^{k-1}(:,\overline{\mathcal J}_1)$ can be extracted efficiently.  The auxiliary sets depend on the type of discretizations, and the notion of dependency of any entry on other entries can be extended to generic TDEs.

Now we present an algorithm for Runge-Kutta time integration for advancing the FOM at DEIM-selected entries.  We consider second-order explicit Runge-Kutta time integration where:
\begin{align*}
   K_1 &=  \mathcal F(\hat{V}^{k-1}),\\
   K_2 &=  \mathcal F(\hat{V}^{k-1} +\frac{\Delta t}{2} K_1),
\end{align*}
and $V^{k} = \hat{V}^{k-1} + \frac{\Delta t}{2}(K_1 + K_2) $. The first TT core requires computing $V^{k}(:,\mathcal J_1) = \hat{V}^{k-1}(:,\mathcal J_1) + \frac{\Delta t}{2}(K_1(:,\mathcal J_1) + K_2(:,\mathcal J_1)) $. While $K_1(:,\mathcal J_1)$ can be computed at the DEIM-selected entries in an identical fashion to the explicit Euler time scheme, the evaluation of  $K_2(:,\mathcal J_1)$ requires providing $\mathcal F(\cdot)$ with $\hat{V}^{k-1}(:,\overline{\mathcal J}_1)+\frac{\Delta t}{2} K_1(:,\overline{\mathcal J}_1)$ as input.  However, $K_1$ is evaluated at $(:,\mathcal J_1)$ entries and the value of $K_1$ at the auxiliary indices is not known.  It is possible to evaluate  $K_1(:,\overline{\mathcal J}_1)$. This approach is unsatisfactory because it requires computing the entries that are  ``auxiliary of the auxiliary".  For higher-order Runge-Kutta schemes the number of auxiliary indices would grow exponentially with the number of Runge-Kutta stages.  To resolve this issue, we build a tensor train low-rank approximation of $K_1$:
\begin{equation}
    \hat{K}_1 = \texttt{TT-CUR-DEIM}(\mathcal F(\hat{V}^{k-1}),\mathcal J_1,\mathcal J_2, \dots, \mathcal J_{d-1}),
\end{equation}
where $\hat{K}_1$ is obtained in TT form $\{G_1,G_2, \dots, G_d \}$ and thus the values of  $\hat{K}_1$ at auxiliary indices can be evaluated efficiently. Note that computing $\hat{K}_1$ utilizes the existing values of $\mathcal F(\hat{V}^{k-1})$ sampled at DEIM-selected entries and it does not require additional evaluation of $\mathcal F(\hat{V}^{k-1})$. Similarly, we construct TT low-rank approximation for the subsequent stages, i.e., $\hat{K}_2, \hat{K}_3$, etc. 
\begin{algorithm}[t]
\begingroup
\fontsize{9pt}{9pt}\selectfont
\caption{\small{\texttt{TT-CUR-DEIM} algorithm for time integration on low-rank tensor train manifold.}}
\label{alg:TT_CUR_DEIM_time}
\SetAlgoLined

\KwIn{$V^0$: function handle that returns the initial condition $V^0(i_1, i_2, \dots, i_d)$ \\

\hspace{10mm} $\mathcal{J}_1, \mathcal{J}_2, \dots,  \mathcal{J}_{d-1}$: randomly initialized right indices.
} \vspace{0.5mm}

\KwOut{$G_1^k, G_2^k, \dots, G_d^k$: core tensors of the tensor train decomposition at $k$th time step}\vspace{2mm}

$[G_1, \dots, G_d;~ \mathcal{I}_1, \dots, \mathcal{I}_{d-1}] =$ \texttt{TT-CUR-DEIM(iterative)}$(V^0, \mathcal{J}_1,  \dots, \mathcal{J}_{d-1})$ \vspace{0.5 mm}

$\mathcal{J}_1^0 \leftarrow \mathcal{I}_{d-1};~~ \mathcal{J}_2^0 \leftarrow \mathcal{I}_{d-2}; ~\dots ~; ~ \mathcal{J}_{d-1}^0 \leftarrow \mathcal{I}_1$ \vspace{0.5 mm}

$G_1^{0} \leftarrow \texttt{permute} \bigl (G_d, [3,2,1] \bigr); \quad \dots \quad; G_d^{0} \leftarrow \texttt{permute} \bigl (G_1, [3,2,1] \bigr)$ \vspace{0.5 mm}

\For{$k = 1: ~$max time-steps}{
[$G_1^k, \dots, G_d^k;~ \mathcal{I}_1^k, \dots, \mathcal{I}_{d-1}^k$] = \texttt{TT-CUR-DEIM}  $\Bigl(\hat{V}^{k-1} + \Delta t \mathcal F(\hat{V}^{k-1}), ~ \mathcal{J}_1^{k-1}, \dots, \mathcal{J}_{d-1}^{k-1} \Bigr)$ \vspace{0.5mm}

$\mathcal{J}_1^{k-1} \leftarrow \mathcal{I}_{d-1}^k;~~ \mathcal{J}_2^{k-1} \leftarrow \mathcal{I}_{d-2}^k; ~\dots ~; ~ \mathcal{J}_{d-1}^{k-1} \leftarrow \mathcal{I}_1^k$ \vspace{0.5 mm}

$G_1^{k-1} \leftarrow \texttt{permute} \bigl (G_d^k, [3,2,1] \bigr); \quad \dots \quad; G_d^{k-1} \leftarrow \texttt{permute} \bigl (G_1^k, [3,2,1] \bigr)$

}
\endgroup
\end{algorithm}

For the time integration of TDE on low-rank manifolds, we can leverage the information from the previous time step to eliminate the iterative process of the \texttt{TT-CUR-DEIM}. This approach is similar to our method for solving MDEs, as explained in \cite{MDOblique}. Specifically, in the time integration of $V^k$, the left-to-right nested indices determined in the previous time step can be reused in the current time step.

In practice, utilizing the indices
from the previous timestep yields excellent performance. The accuracy difference when \texttt{TT-CUR-DEIM} is used iteratively is negligible. A similar approach was employed in \cite{MNAdaptive, MDOblique, GB24}, where information from the previous timestep is utilized in the current timestep to solve MDEs and TDEs. However, at $t=0$, \texttt{TT-CUR-DEIM (iterative)} must be used to determine the left or right-nested indices.  The time integration algorithm for solving TDEs in tensor train low-rank form is presented in Algorithm \ref{alg:TT_CUR_DEIM_time}.

\subsection{Rank adaptivity}
\label{Rank adaptivity}\label{sec:rank-adaptivity}
The proposed time integration algorithm can be augmented with rank adaptivity, enabling on-the-fly control of the error in the low-rank approximation. To this end, we define an error proxy as:
\begin{equation} \label{ErrorCriterion}
    \epsilon_z = \dfrac{\min (\boldsymbol \Sigma_z)}{\|\boldsymbol \Sigma_z \|_F}  \quad z=1,\dots,d-1,
\end{equation}
where $\boldsymbol \Sigma_z$ is the matrix of singular values that can be obtained in the first step of Algorithm \ref{alg:CUR-DEIM}  as $[\bm U , \boldsymbol \Sigma_z , \sim ] = \texttt{SVD}(V(:,\mathcal J), r_z), ~ z=1,\dots, d-1$, where $\mathcal J$ is the appropriate set of indices. The value of $\epsilon_z$  measures the relative contribution of the $r$-th rank. The rank $r_z$ is adjusted or remains unchanged to maintain $\epsilon$ within a desired range, $\epsilon_{l} \leq \epsilon_z \leq \epsilon_{u}$, where $\epsilon_{l}$ and $\epsilon_{u}$ are user-specified lower and upper thresholds. In case $\epsilon_z < \epsilon_l$, the rank is reduced by one such that the new rank is $r'_z = r_z -1$. Once the new rank is determined, the DEIM-selected indices $\mathcal{I}_z$ are truncated to retain the first $r'_z$ components. On the other hand, when $\epsilon_z > \epsilon_u$, the new rank is obtained as $r'_z = r_z+1$. The rank is increased by sampling more indices and updating $\mathcal{I}_z$. Since DEIM only provides sampling points equal to the number of the columns of the left singular matrix (i.e. $\bm U$), we perform oversampling using GappyPOD+E algorithm presented in \cite{BPStability}. Although any sparse selection method can be used, it has been demonstrated that GappyPOD+E outperforms alternative approaches such as random sampling or leverage scores \cite{MMCUR}. By utilizing the GappyPOD+E, we add $m$ more indices to the index sets $\mathcal{I}_1, \mathcal{I}_2, \dots \mathcal{I}_{d-1}$, where $m \geq 1$ is an arbitrary number. As $m$ increases, the low-rank approximation error will decrease, however, in practice, $m=5$ can capture the desired level of accuracy. Although the algorithm adds $m$ more indices to the index sets, the rank is increased by one since we truncate the SVD of the function evaluations at $r+1$ (see the first step of Algorithm \ref{alg:CUR-DEIM}).

\subsection{Extension to  Tucker tensor low-rank approximation}
In this section, we present an algorithm for solving TDE on low-rank Tucker tensor manifolds. While the primary focus of this paper is on tensor train low-rank approximation, the extension to the Tucker tensor model, which is distinctly different from tensor train low-rank approximation, is another demonstration of the utility of the cross algorithms for the time-integration of TDEs on low-rank manifolds. 

We demonstrate that Eq. \ref{TDE-CROSS} can result in a cost-effective and well-conditioned time integrator for the Tucker tensor model provided that an effective cross algorithm is utilized. To this end, we use a Tucker tensor cross algorithm based on DEIM that we recently developed \cite{GB24}. This algorithm is referred to as DEIM fiber sampling or \texttt{DEIM-FS}. In \cite{GB24}, the \texttt{DEIM-FS} algorithm was utilized to reduce the computational cost of solving DLRA evolution equations, which are based on tangent space projection and can become unstable in the presence of small singular values. In this work, we utilize the \texttt{DEIM-FS} algorithm to judiciously sample the time-discrete FOM at the DEIM-selected fibers, which as we show, results in a cost-effective and stable time integration scheme. Another excellent candidate for the Tucker cross algorithm is presented in \cite{DKS21}, which could be utilized for the time integration of TDEs on low-rank Tucker tensor maniflds. 

The \texttt{DEIM-FS} algorithm is presented in Algorithm \ref{alg:DEIM-FS} for a third-order tensor. We briefly review this algorithm here and refer the reader to \cite{GB24} for more details. The input to the \texttt{DEIM-FS} algorithm is the approximate or exact singular vectors, denoted by $\tilde{\bm U}_m$, of mode$-m$ unfolded matrices, i.e., $\bm V_{(m)}$. In Step 1, DEIM is applied to $\tilde{\bm U}_m$, which yields the DEIM indices along each mode, i.e., $\mathcal I_m$. In the case of solving TDEs, the singular vectors of the previous time step are used as $\tilde{\bm U}_m$.   In Step 2, the target tensor is evaluated at the DEIM-selected fibers. In the case of solving TDEs, the time-discrete FOM (Eq. \ref{TimeDiscreteTDE}) is sampled at these fibers. For Runge-Kutta time integration, this step involves building the low-rank Tucker model for $K_i$ tensors, which can then be used for computing the values of  $K_i$ at auxiliary indices. The sampled fibers along mode $m$ form the $\bm C_m$ matrices. The SVD of these matrices and truncation at rank $r_m$ produces the factorized matrices, $\bm U_m$. In Step 4, the intersection tensor ($W$) is formed. The core tensor ($G$) is computed by interpolating the target tensor onto the space spanned by the columns of $\bm U_m$ at the DEIM intersection indices, which is shown in Step 5. In summary,  the \texttt{DEIM-FS} algorithm produces $\{G,\bm U_1, \bm U_2, \bm U_3\}$, which are the components of the low-rank Tucker model $\hat{V}$ according to Eq. \ref{eq:TUform}. As shown in \cite[Remark 1]{GB24}, $\hat{V}(\mathcal I_1, \mathcal I_2, \mathcal I_3) = V(\mathcal I_1, \mathcal I_2, \mathcal I_3)$. It is possible to oversample $V$ at $r'_m > r_m$ fibers which can result in more accurate approximations.  
 
The above steps are in several ways analogous to the \texttt{TT-CUR-DEIM} algorithm as well as to the CUR approach presented in \cite{MDOblique}. In all of these methods, DEIM indices are determined using the singular vectors from the previous time step, and low-rank approximations are constructed for Runge-Kutta $K_i$ tensors/matrices. Additionally, the initial condition can be computed in Tucker tensor low-rank format using the \texttt{DEIM-FS (iterative)} algorithm from \cite{GB24}, which is analogous to \texttt{TT-CUR-DEIM (iterative)} and it does not require access to the singular vectors.

\begin{algorithm}[t]
\begingroup
\fontsize{9pt}{9pt}\selectfont
\caption{\small{\texttt{DEIM-FS} algorithm for 3D Tucker tensor low-rank approximation \cite{GB24}.}}
\label{alg:DEIM-FS}
\SetAlgoLined

\KwIn{$\tilde{\bm U}_m$:  matrix of exact or approximate left singular vectors of $V_{(m)}$\\

\hspace{11mm} $r_{m}$:  target Tucker rank\\

\hspace{11mm} $V$:  function handle to compute fibers of the target $V$ }

\KwOut{$G,~ \bm U_1 ,~ \bm U_2,~ \bm U_3$}\vspace{3mm}

$\mathcal{I}_m=$ \texttt{DEIM}($\tilde{\bm U}_m$)  \hspace{85mm} $\rhd$ compute $\mathcal{I}_1$, $\mathcal{I}_2$, $\mathcal{I}_3$ \vspace{1mm}

$\bm C_1 = \Bigl(V(:, \mathcal{I}_2, \mathcal{I}_3)\Bigr)_{(1)}$ \hspace{28mm} $\rhd$ calculate $\bm C_1 \in\mathbb{R}^{n_1\times r_{2}r_{3}}$, ~$\bm C_2 \in\mathbb{R}^{n_2\times r_{1}r_{3}}$, $\bm C_3 \in\mathbb{R}^{n_3\times r_{1}r_{2}}$ \vspace{-1mm}

\nonl $\bm C_2 = \Bigl(V(\mathcal{I}_1, :, \mathcal{I}_3)\Bigr)_{(2)}$ \vspace{-1mm}

\nonl $\bm C_3 = \Bigl(V(\mathcal{I}_1, \mathcal{I}_2, :)\Bigr)_{(3)}$
\vspace{2mm}

$[\bm U_m, \sim , \sim ] = \texttt{SVD}( \bm C_m, r_{m})$ \hspace{8mm}  $\rhd$ calculate  the left singular vectors of $\bm C_i$ and truncate at rank $r_{m}$ \vspace{2mm}

$W = V(\mathcal{I}_1, \mathcal{I}_2, \mathcal{I}_3)$
\hspace{50mm} $\rhd$ form the intersection tensor $W \in\mathbb{R}^{r_{1}\times r_{2}\times r_{3}}$ \vspace{2mm}

$G = W \times_1 ~ \bm U_1(\mathcal{I}_1, :)^\dagger \times_2 ~ \bm U_2(\mathcal{I}_2, :)^\dagger \times_3 ~ \bm U_3(\mathcal{I}_3, :)^\dagger$ \hspace{26mm} $\rhd$ calculate the core tensor ($G$) \vspace{1mm}

\endgroup
\end{algorithm}

\subsection{Summary of computational advantages}
The above algorithm is remarkably simple and can be explained in the following steps: The tensor $V^{k} = \hat{V}^{k-1} + \Delta t \overline{\mathcal F}(\hat{V}^{k-1})$ represents the solution of the FOM at time step $k$, with the solution at the previous time step, $k-1$, set as the low-rank tensor. However, $V^k$, which would be high or full-rank, is never actually formed. Instead, the low-rank tensor at time step $k$ is obtained by performing a rank truncation using a cross algorithm, which judiciously evaluates $V^{k}$ at sparse entries.

 In the following, we list the advantages of the above time-integration algorithm:

(i) \textbf{Computational Cost:} For tensor train decomposition, \texttt{TT-CUR-DEIM} only samples $dnr^2$ entries of the cross algorithm regardless of the type of nonlinearity of $\mathcal F(\ \cdot \ )$. This means that the  \texttt{TT-CUR-DEIM} evaluates the right-hand side function at the minimum number of entries possible for a rank-$r$ approximation, which is significantly smaller than the full-size tensor, which requires $n^d$ function evaluations. The algorithm is also memory efficient and has the storage complexity of $\mathcal{O}(dnr^2)$.  The \texttt{DEIM-FS} algorithm  has memory and cost complexity of $\mathcal{O}(r^d+ndr^{d-1})$. 

(ii) \textbf{Ill-Conditioning:} The inversion of the singular value matrix is not needed in the above algorithm and it is robust in the presence of small or zero singular values. Similar to what was shown in \cite{MDOblique}, the tensor DEIM-based cross algorithms are well-conditioned. The same is true for maxvol cross algorithms  \cite{GTZ97,ISTtcross}.

(iii) \textbf{Runge-Kutta Temporal Integration:} The presented algorithm inherits the temporal accuracy of the scheme used for temporal integration. Details of the Runge-Kutta time integration scheme are presented in Section \ref{high-orderTime}. As shown in Figures \ref{fig:DynamicToyTT} and \ref{fig:DynamicToyTucker}, the time integration method presented for both tensor train and Tucker tensor decompositions retains the fourth-order temporal accuracy of RK4 for nonlinear TDEs.

(iv) \textbf{Intrusiveness:} The above algorithm is easy to implement, mainly because it only evaluates the full-order model at judiciously chosen entries. It does not involve tangent space projections, nor does it require term-by-term treatment of the low-rank tensor, both of which can be quite involved depending on the TDE. The cross algorithms are implemented only once and are straightforward to implement. For the most part, they are agnostic to the type of TDEs being applied. All the examples in this paper utilize the same cross algorithm.

(v) \textbf{Rank Adaptivity:} It is straightforward to perform rank adaptivity using the cross algorithms, which involves sampling additional or fewer fibers if rank addition or removal criterion is triggered.  We present the details of rank adaptivity in Section  \ref{Rank adaptivity}.

\section{Demonstrations}

\subsection{Toy Examples} 
\label{sec:ToyEx}

As the first example, we consider two three-dimensional functions as shown below:
\begin{subequations}
\begin{align}
  \mathcal F_1(x_1,x_2,x_3)  = e^{-4(x_1 x_2 x_3)^2} \hspace{3.4cm}  x_1, x_2, x_3 \in [-1,1] ,
  \label{eq:toy1} \\
  \mathcal F_2(x_1,x_2,x_3)  = \frac{1} {\bigl(x_1^b + x_2^b + x_3^b\bigr)^{1/b}} \hspace{1cm} x_1, x_3 \in [1,200]  ~~~~ x_2 \in[1,300] .
  \label{eq:toy2}
\end{align}
\end{subequations}

\begin{figure}[ht]
     \centering
     \begin{subfigure}[t]{0.45\textwidth}
         \centering
         \includegraphics[scale=1]{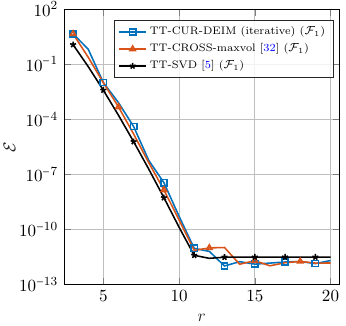}
         \caption{}
         \label{fig:F1Errors}
     \end{subfigure}
     \begin{subfigure}[t]{0.45\textwidth}
         \centering
         \includegraphics[scale=1]{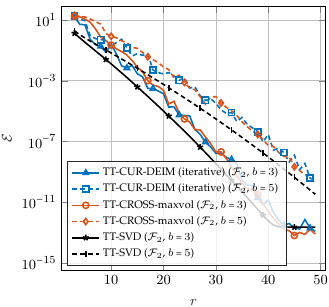}
         \caption{}
         \label{fig:F2Errors}
     \end{subfigure}
     \\
    \begin{subfigure}[l]{0.45\textwidth}
         \centering
         \includegraphics[scale=1]{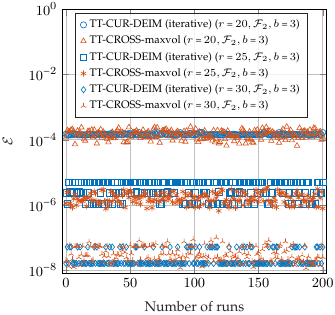}
         \caption{}
         \label{fig:randomInitial}
     \end{subfigure}
     \caption{Toy Examples: Comparison of the \texttt{TT-CUR-DEIM (iterative)}, TT-Cross-maxvol, and TT-SVD: (a) Approximation error of $F_1$ versus rank; (b) Approximation error of $F_2,  (b=3,5)$ versus rank; (c) Effect of random fiber initialization on the approximation error of $F_2, (b=3)$ for three different ranks ($r=20$, $r=25$, $r=30$).}
     \label{fig:CompareMethods}
\end{figure}

By evaluating $\mathcal{F}_1(x_1, x_2, x_3)$ at 100 equally spaced points for each of $x_1$, $x_2$, and $x_3$ within their respective domains, we obtain $F_1$, a tensor of size $100 \times 100 \times 100$. Similarly, evaluating $\mathcal{F}_2(x_1, x_2, x_3)$ on an equally spaced grid in each direction results in a tensor $F_2 \in \mathbb{R}^{200 \times 300 \times 200}$. Two choices of $b=3$ and $b=5$ are considered for $F_2$. To compare the performance of the DEIM and maxvol algorithms, $F_1$ and $F_2$ are approximated using \texttt{TT-CUR-DEIM} and TT-Cross-maxvol, respectively. Then, the error between the approximated tensors and the actual tensors is calculated. The low-rank approximation error is defined as:
\begin{equation*}
 \mathcal{E}  = \| \hat{F} - F \|_F ,
\end{equation*}
where $\hat{F}$ denotes the low-rank tensor train approximation of $F$. For both methods, we consider $r_1 = r_2 = r_3 =r$. Figures \ref{fig:F1Errors} and \ref{fig:F2Errors} show the low-rank approximation error versus rank ($r$) for $F_1$ and $F_2$, respectively. As shown in these figures, the error of TT-Cross-maxvol algorithm closely follows the \texttt{TT-CUR-DEIM} error. The results demonstrate that both DEIM and maxvol methods are performing comparably and the obtained errors are close to the error of TT-SVD \cite{TT11}.

In both \texttt{TT-CUR-DEIM} and TT-Cross-maxvol methods, the initial right-nested indices are selected randomly. We conducted a study to examine the effect of different initialization on the accuracy of the low-rank approximations. Figure \ref{fig:randomInitial} depicts the approximation error of $F_2$ ($b=3$) obtained from 200 random initializations of both \texttt{TT-CUR-DEIM} and TT-Cross-maxvol algorithms. Overall both techniques show similar errors but DEIM shows a smaller variance in most cases considered. 

Figure \ref{fig:randomInitial} illustrates that depending on the rank and initialization, one of the algorithms may result in a slightly lower error compared to the other. However, in general, they both exhibit similar accuracy.

\subsection{High-dimensional nonlinear TDE}
\label{sec:Example2}

\begin{figure}[t]
     \centering
     \begin{subfigure}[t]{0.45\textwidth}
         \centering
         \includegraphics[scale=1]{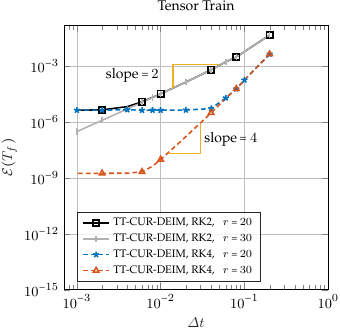}
         \caption{}
         \label{fig:DynamicToyTT}
     \end{subfigure}
     \begin{subfigure}[t]{0.45\textwidth}
         \centering
         \includegraphics[scale=1]{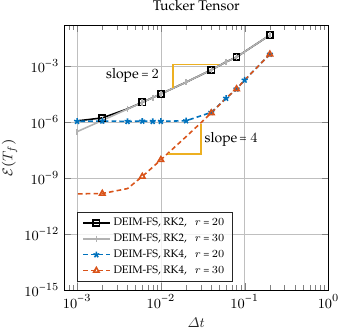}
         \caption{}
         \label{fig:DynamicToyTucker}
     \end{subfigure}
     \\
     \begin{subfigure}[t]{0.45\textwidth}
         \centering
         \includegraphics[scale=1]{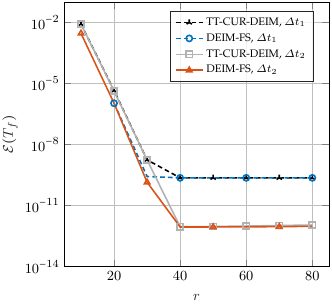}
         \caption{}
         \label{fig:ToyTDE_rConvereges}
     \end{subfigure}
     \caption{Three-dimensional nonlinear TDE: Error at the final time ($\mathcal{E}(T_f)$) in solving the given TDE in Eq. \ref{TDE toy example} versus time step size ($\Delta t$). The second-order and fourth-order Runge-Kutta (RK2 and RK4) are used for time integration with various ranks $r=20$, $r=30$. (a) \texttt{TT-CUR-DEIM}; (b) \texttt{DEIM-FS}; (c) Error at the final time ($\mathcal{E}(T_f)$) versus rank with $\Delta t_1 = 4\times 10^{-3}$ and $\Delta t_2 = 1\times 10^{-3}$.}
     \label{fig:CompareTTandTucker}
\end{figure}
One of the important advantages of applying residual minimization at the time discrete level is that the resulting low-rank approximation inherits the order of accuracy of the temporal discretization of the FOM. To illustrate this aspect of the methodology, we consider a time-dependent tensor given by: 
\begin{equation} \label{dynamic toy example}
    v(x_1, x_2, \dots, x_d, t) = \frac{1} {\bigl(\sum_{i=1}^d x_i^b + e^{\lambda t}\bigr)^{1/b}} \hspace{1cm} x_i \in [1,200], \quad i=1,2, \dots, d,  \quad t \in[0,1] .
\end{equation}
We first consider $d=3$ where the time-dependant tensor $V (t) \in \mathbb{R}^{200 \times 200 \times 200}$ is obtained by evaluating $v(x_1, x_2, x_3, t)$ on a uniform grid of $x_1$, $x_2$, $x_3$ at a particular time (t). Based on $V(t)$, a nonlinear TDE can be obtained as below:
\begin{equation} \label{TDE toy example}
    \frac{d V}{dt} = -\frac{\lambda ~ e^{\lambda t}}{b} V^{b+1}.
\end{equation}
We use explicit second-order and fourth-order Runge-Kutta integrators (RK2 and RK4) for temporal integration of the above TDE. We also considered both Tucker tensor and tensor train low-rank approximation schemes using \texttt{TT-CUR-DEIM} and \texttt{DEIM-FS} algorithms, respectively. We choose $b=3$, $\lambda=10$ and final time $T_f = 1$. We use the relative error in this example and the rest of the subsequent examples, which is defined as:
\begin{equation}\label{eq:err_rel}
 \mathcal{E}(t)  = \frac{\| \hat{V}(t) - V(t) \|_F~}{~\| V(t) \|_F}.
 \end{equation}
 In Figures \ref{fig:DynamicToyTT} and \ref{fig:DynamicToyTucker}, the obtained error at the final time ($\mathcal{E}(T_f)$) is depicted versus time step size ($\Delta t$) for \texttt{TT-CUR-DEIM} and \texttt{DEIM-FS} models. Two fixed-rank values of $r=20, 30$ are considered. The figures illustrate that both \texttt{TT-CUR-DEIM} and \texttt{DEIM-FS} methods retain the fourth-order accuracy characteristic of the RK4 scheme, as well as the second-order accuracy of the RK2 scheme. As anticipated, the \texttt{TT-CUR-DEIM} and \texttt{DEIM-FS} with RK4 saturates to the low-rank error for each $r$ faster than RK2. This example demonstrates that the DEIM cross algorithms inherit the temporal order of accuracy.

 In Figure \ref{fig:ToyTDE_rConvereges}, we plot the error at the final time ($\mathcal{E}(T_f)$) versus rank for two different step sizes, $\Delta t_1 = 4 \times 10^{-3}$ and $\Delta t_2 = 1 \times 10^{-3}$. We consider the \texttt{TT-CUR-DEIM} and \texttt{DEIM-FS} methods for solving the TDE on the low-rank tensor train and Tucker tensor manifolds. According to the figure, $\mathcal{E}(T_f)$ decreases as the rank increases until the temporal integration error dominates. When the time step is reduced to $\Delta t_2$, we observe a corresponding decrease in the saturated errors. As the rank increases, the smallest singular values of the unfolded matrices decrease to near machine precision. Figure \ref{fig:ToyTDE_rConvereges} verifies that both the \texttt{TT-CUR-DEIM} and \texttt{DEIM-FS} methods are robust against the presence of small singular values, even when the rank is unnecessarily large, i.e., overapproximation.

\begin{figure}[htbp]
     \centering
     \begin{subfigure}[t]{0.45\textwidth}
         \centering
         \includegraphics[scale=1]{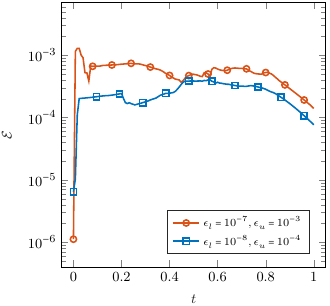}
         \caption{}
         \label{fig:High_dim_error}
     \end{subfigure}
     \begin{subfigure}[t]{0.45\textwidth}
         \centering
         \includegraphics[scale=1]{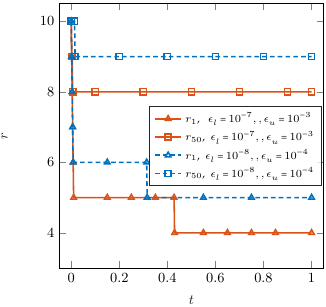}
         \caption{}
         \label{fig:High_dim_ranks}
     \end{subfigure}
     \caption{100-dimensional nonlinear TDE on low-rank tensor train manifold: (a) Relative error evolution; (b) Evolution of $r_1$ and $r_{50}$ (the first and 50-th tensor train ranks). }
     \label{fig:100D_error}
\end{figure}

To demonstrate the performance of the proposed \texttt{TT-CUR-DEIM} method in dealing with high-dimensional tensor differential equations, we consider a 100-dimensional case ($d=100$) of the given TDE in Eq. \ref{TDE toy example}. We consider the time-dependant tensor $V (t) \in \mathbb{R}^{70 \times 70 \times \dots \times 70}$ is obtained by evaluating $v(x_1, \dots, x_{100}, t)$ on a uniform grid of $x_1, \dots, x_{100}$. Therefore, the total number of tensor elements is $70^{100} \approx 3.2 \times 10^{184}$. We use RK4 for temporal integration using \texttt{TT-CUR-DEIM}. We choose b = 0.9, which corresponds to a non-polynomial nonlinearity and $\Delta t = 2\times 10^{-3}$. The TDE is solved with two different values of $\epsilon_l$ and $\epsilon_u$. In this example and the following examples, the difference between $\epsilon_l$ and $\epsilon_u$ must be large enough to avoid oscillation in rank adaptivity.  Because of the size of the tensor, we could not calculate the Frobenius norm that requires a summation over all entries. Instead, we
calculate the error at roughly 3,000 random entries.  Figure \ref{fig:High_dim_error} shows the relative error over time for $\epsilon_l = 10^{-7}, \epsilon_u = 10^{-3} $ and $\epsilon_l = 10^{-8}, \epsilon_u = 10^{-4}$. Figure \ref{fig:High_dim_ranks} depicts the rank versus time based on the selected $\epsilon_l$ for solving the TDE. Due to the non-polynomial nonlinearity, the existing techniques for solving TDE on the low-rank manifolds are impractical for this problem due to memory and FLOP costs. However, \texttt{TT-CUR-DEIM} enables compressing and integrating a tensor of size  $3.2\times 10^{184}$ with effectively $dnr^2 \approx 100\times 70 \times 10^2 = 7 \times 10^5$ entries.

\subsection{Four-dimensional nonlinear advection-reaction (4D-AR) equation}
\label{sec:Example3}
As the next demonstration, we consider a four-dimensional nonlinear advection-reaction  (4D-AR) given below:
\begin{equation}\label{eq:NonlinearPDE}
 \frac{\partial v(\bm x,t)}{\partial t} = -\bm b \cdot \nabla v(\bm x,t) + s(v(\bm x,t)), \quad \bm x \in [-5,5]^4,
\end{equation} 
where $\bm b = [-\sin(t),~ \cos(t),~ -\sin(\pi + t),~ \cos(\pi + t)]$ is the advection velocity and $s(v)=- 0.1e^{-v}/(1+v^2)$ is the nonlinear source term.  These types of nonlinearity are encountered for example in DLRA of turbulent combustion \cite{RNB21}. The boundary condition is periodic in all directions and the initial condition is considered as: 
\begin{equation*}
 v_0(\bm x) = e^{-(2x_1-\frac{1}{2})^2}  e^{-(2x_2-\frac{1}{2})^2} e^{-(2x_3-\frac{1}{2})^2} e^{-(2x_4-\frac{1}{2})^2}.
\end{equation*}
The 4D-AR PDE, given in Eq. \ref{eq:NonlinearPDE}, includes non-polynomial nonlinear terms on the right-hand side. In this case, $\mathcal{F}(\hat{V})$ remains a full-rank tensor, despite $\hat{V}$ being of low-rank. Consequently, the implementation of the standard DLRA necessitates the explicit formation of the full-dimensional $\mathcal{F}(\hat{V})$. This would demand substantial memory allocation and FLOP costs that scale with $\mathcal{O}(n^4)$. The \texttt{TT-CUR-DEIM} and \texttt{DEIM-FS} algorithms circumvent this cost, as the full-rank tensor is not formed.

We employ a spectral-element method for the spatial discretization ($\bm x$). Each spatial domain is discretized using $n=n_1 = n_2 = n_3 = n_4 = 81$ points, consisting of 40 elements with a second-order polynomial approximation within each element. The RK4 method is used for temporal integration within the time interval $t \in [0,4]$ with $\Delta t = 1 \times 10^{-3}$. The \texttt{TT-CUR-DEIM} and \texttt{DEIM-FS} algorithms are employed to solve the TDE in the low-rank tensor train and Tucker tensor forms. For error analysis, the relative error as defined by Eq. \ref{eq:err_rel} is used. As the analytical solution is unavailable for this TDE, the solution obtained from the FOM is considered the ground truth. The FOM is solved with the same spatiotemporal resolution as the low-rank model.

\begin{figure}[t]
     \centering
     \begin{subfigure}[t]{0.45\textwidth}
         \centering
         \includegraphics[scale=1]{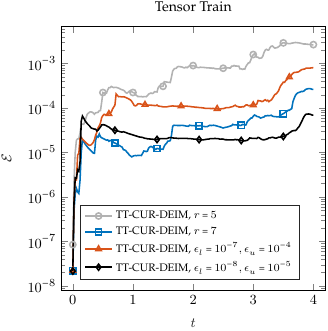}
         \caption{}
         \label{fig:Ex2Error}
     \end{subfigure}
     \begin{subfigure}[t]{0.45\textwidth}
         \centering
         \includegraphics[scale=1]{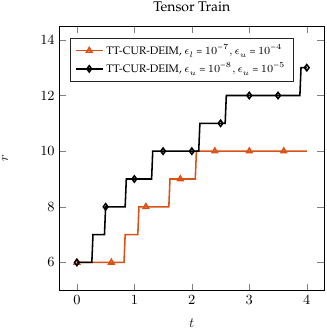}
         \caption{}
         \label{fig:Ex2rankevolution}
     \end{subfigure}
     \\
     \begin{subfigure}[t]{0.45\textwidth}
         \centering
         \includegraphics[scale=1]{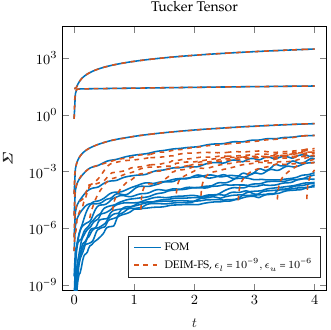}
         \caption{}
         \label{fig:TuAdvecionSingval}
     \end{subfigure}
     \begin{subfigure}[t]{0.45\textwidth}
         \centering
         \includegraphics[scale=1]{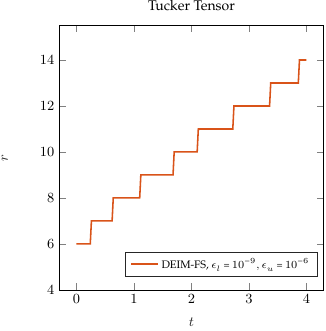}
         \caption{}
         \label{fig:TuAdvecionrank}
     \end{subfigure}
     \caption{Four-dimensional nonlinear advection equation on low-rank tensor manifolds: (a)  Relative error evolution with fixed rank and adaptive rank solved with \texttt{TT-CUR-DEIM}; (b) Evolution of $r$ associated with the rank-adaptive systems plotted in Figure \ref{fig:Ex2Error}; (c) Singular values of the mode-1 unfolding matrix obtained from FOM and \texttt{DEIM-FS}; (d) Evolution of $r=r_1=r_2=r_3=r_4$ associated with the rank-adaptive case plotted in Figure \ref{fig:TuAdvecionSingval}. }
     \label{fig:CompareErros}
\end{figure}

We first solve the TDE using the \texttt{TT-CUR-DEIM} algorithm with the fixed rank of $r=5$ and $r=7$. Figure \ref{fig:Ex2Error} indicates that as the rank increases from $r=5$ to $r=7$, the error decreases. We also solve the TDE with two rank-adaptivity thresholds of $\epsilon_l = 10^{-7}, \epsilon_u = 10^{-4}$ and $\epsilon_l = 10^{-8}, \epsilon_u = 10^{-5}$. 

As shown in Figure \ref{fig:Ex2rankevolution}, a lower value of the upper threshold ($\epsilon_u$) results in larger ranks and a more accurate approximation.

The same TDE is also solved on Tucker tensor low-rank manifolds using the \texttt{DEIM-FS} algorithm with $\epsilon_l = 10^{-9}, \epsilon_u = 10^{-6}$. Figure \ref{fig:TuAdvecionSingval} illustrates a comparison between the singular values of the mode-1 unfolding, i.e., $\hat{\bm V}_{(1)}$ obtained from the \texttt{DEIM-FS} algorithm and those of FOM, while Figure \ref{fig:TuAdvecionrank} presents the corresponding evolution of $r=r_1=r_2=r_3=r_4$ for this problem. 

\subsection{Stochastic advection-diffusion-reaction (S-ADR) equation} 
\label{sec:3D_ADR}
In this section, we consider a three-dimensional nonlinear stochastic advection-diffusion-reaction equation subject to random diffusion and reaction coefficients $\xi_1$ and $\xi_2$, with deterministic initial condition:  
\begin{equation}\label{eq:3D_ADR}
 \frac{\partial v}{\partial t} + (u \cdot \nabla ) v = \nabla \cdot (\nabla v \xi_1) + \frac{v^2}{10 + \xi_2 v}, \quad x_1 \in [0,10],~ x_2 \in [0,2],~ t \in [0,5],
\end{equation} 
\begin{equation*}
v(\bm x,0; \xi_1,\xi_2) = \frac{1}{2} \biggl ( \tanh \Bigl (\frac{x_2 + 0.5}{0.1} \Bigr) - \tanh \Bigl(\frac{x_2 - 0.5}{0.1} \Bigr) \biggr),
\end{equation*} 
where $\bm x=(x_1,x_2)$, $v(\bm x, t; \xi_1, \xi_2)$ is the species concentration and $u(\bm x, t)$ is the deterministic velocity vector. This example also includes non-polynomial nonlinear terms in the RHS. The schematic of the problem is presented in our previous studies \cite{MDOblique, MDComputing, RNB21}. The velocity field is obtained by solving the incompressible Navier–Stokes equations and is independent of the species transport equation.  In particular, we solve the velocity field in the entire domain using the spectral/hp element method. At the inlet, a parabolic velocity is prescribed, with an average velocity of $\overline{u}$. The outflow condition is imposed at the right boundary and the no-slip boundary condition is imposed at the remaining boundaries. The Reynolds number with reference length $H/2$, and kinematic viscosity $\nu$, is given by $Re = \overline{u} H/2\nu = 1000$.

We formulate the above PDE as a three-dimensional TDE, where the TT low-rank representation in the continuous form is given by:
\begin{equation*}
    v(\bm x,t;\xi_1,\xi_2) \approx \sum_{\alpha_1=1}^{r_1}\sum_{\alpha_2=1}^{r_2}  g_1(\bm x,\alpha_1;t) g_2(\alpha_1,\xi_1,\alpha_2;t)g_3(\alpha_2,\xi_2;t),
\end{equation*}
where $g_k$'s represent the time-dependent tensor cores in the functional form.
Therefore, the first core contains the two-dimensional spatial functions, and in discrete form, this amounts to vectorizing the two-dimensional spatial grid into one long index.

For the spatial discretization, we use a uniform quadrilateral mesh with 50 elements in the \(x_1\) direction and 15 elements in the \(x_2\) direction, and a spectral polynomial of order 5 in each direction. This results in \(n_1 = 19076\) degrees of freedom along the first mode of the tensor. The fourth-order explicit Runge-Kutta method is used for time integration with \(\Delta t = 5 \times 10^{-4}\).

\begin{figure}[t]
     \centering
     \begin{subfigure}[t]{0.45\textwidth}
         \centering
         \includegraphics[scale=1]{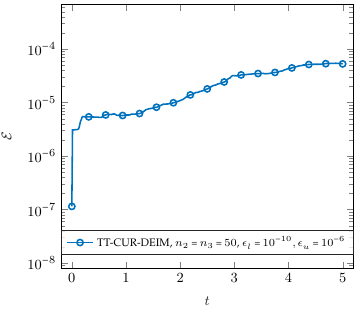}
         \caption{}
         \label{fig:3D_ADR_error}
     \end{subfigure}
     \begin{subfigure}[t]{0.45\textwidth}
         \centering
         \includegraphics[scale=1]{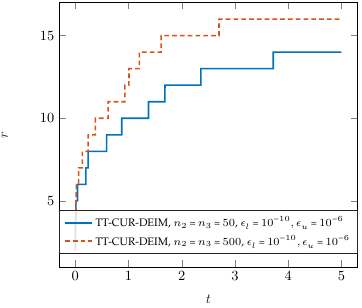}
         \caption{}
         \label{fig:3D_ADR_rank}
     \end{subfigure}
     \caption{Stochastic advection-diffusion-reaction equation on low-rank tensor train manifolds: (a)  Relative error evolution solved with \texttt{TT-CUR-DEIM}  and $n_2 =n_3= 50$; (b) Evolution of $r$ associated with the rank-adaptive systems with $n_2 =n_3= 50$ and $n_2 =n_3= 500$. In both cases, $n_1= 19076$.}
     \label{fig:3D_ADR}
\end{figure}

For the first case, we consider a random diffusion coefficient defined as \(\xi_1 = 1/\alpha_1\), where \(\beta_1\) is a Gaussian random variable with a mean of 100 and a standard deviation of 25, i.e.,  $\beta_1 \sim \mathcal{N}(100, 25)$. The second random coefficient, \(\xi_2\), is parameterized as \(1/\alpha_2\), where \(\beta_2\) is randomly selected from a Gaussian distribution  of $\beta_2 \sim \mathcal{N}(1, 0.05)$. To compute the error of the TT-CUR-DEIM method, we solve the FOM as the ground truth. Since considering a large number of samples (\(n_2, n_3\)) makes the FOM too expensive to compute, we initially consider \(n_2 = n_3 = 50\) samples. The solution tensor is of size \(V(t) \in \mathbb{R}^{n_1 \times n_2 \times n_3}\). Figure \ref{fig:3D_ADR_error} presents the calculated relative error between the FOM and TT-CUR-DEIM with \(\epsilon_l = 10^{-10}, \epsilon_u = 10^{-6}\). The corresponding rank evolution plot is depicted in Figure \ref{fig:3D_ADR_rank}.

\begin{figure}[t]
 \centering
\includegraphics[width=\textwidth]{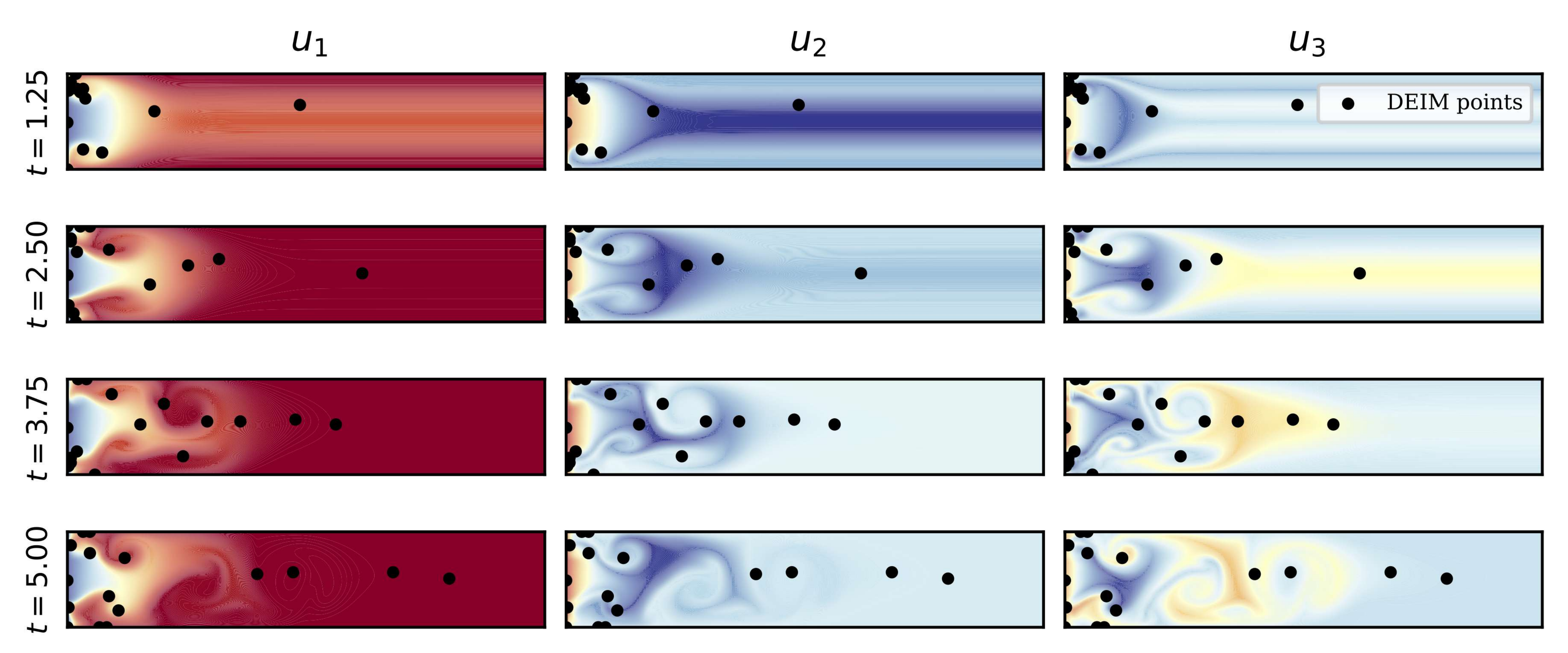}
 \caption{Three-dimensional stochastic advection-diffusion-reaction equation on low-rank tensor train manifolds: First three spatial modes at different time-steps and the selected points of DEIM with $n_2 =n_3= 500$.}
 \label{fig:3D_ADR_solution}
\end{figure}

In this example, we use a deterministic initial condition. Therefore, the rank at $t=0$ is exactly equal to one. While it is possible to initialize the \texttt{TT-CUR-DEIM} with $r=1$, however, we intentionally initialize the simulation with $r=2$ at $t=0$. With this initialization rank, one of the singular values is exactly equal to zero. Despite having zero singular values at $t=0$, the \texttt{TT-CUR-DEIM} remains stable.

For the second case of S-ADR, we take $n_2=n_3 = 500$ samples. In this case, FOM is prohibitively costly to solve. The rank versus time plot for this simulation is shown in Figure \ref{fig:3D_ADR_rank}. The full-order tensor for this case is of size $V(t) \in \mathbb{R}^{19076\times 500 \times 500}$, which is roughly equal to $4.7\times 10^9$ entries.  In Figure \ref{fig:3D_ADR_solution}, we present the evolution of the first three spatial modes, along with the DEIM-selected points. These first three spatial modes are obtained from the first three columns of $\bm U_1$. The figure indicates that as the simulation evolves in time, the points also evolve as the flow is advected from left to right. 

\subsection{Computational complexity} 
In this section, we present the results of the computational complexity analysis.  
All the simulations presented in this paper are executed on a computer with 13-th Gen Intel(R) Core(TM) i9-13900K 3.00 GHz CPU and 32.0 GB RAM. 

\begin{figure}[t]
     \centering
     \begin{subfigure}[t]{0.45\textwidth}
         \centering
         \includegraphics[scale=1]{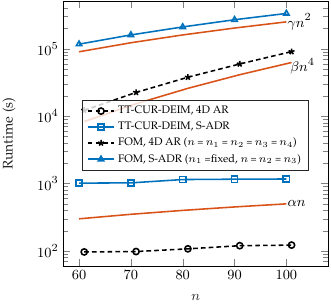}
         \caption{}
         \label{fig:ComputCost}
     \end{subfigure}
     \begin{subfigure}[t]{0.45\textwidth}
         \centering
         \includegraphics[scale=1]{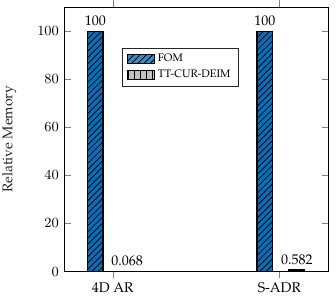}
         \caption{}
         \label{fig:Barchart}
     \end{subfigure}
     \caption{Four-dimensional nonlinear advection-reaction equation (4D-AR) and stochastic advection-diffusion-reaction equation (S-ADR) on low-rank tensor train manifolds: (a) Computational times vs $n$ (number of elements in different tensor modes); (b) Comparison of the relative memory requirements for solving the 4D-AR and S-ADR TDEs.}
     \label{fig:SpeedUps}
\end{figure}

Figure \ref{fig:ComputCost} illustrates the computational cost of solving 4D-AR TDE and S-ADR TDE  versus $n$ for FOM and \texttt{TT-CUR-DEIM}. In the S-ADR example, we only vary $n=n_2=n_3$ while $n_1=19076$ is kept constant.  Therefore, for S-ADR, the cost of solving FOM scales quadratically as $ n$ increases while the computational cost of \texttt{TT-CUR-DEIM} increases linearly.  For 4D-AR, as $n=n_1=n_2=n_3=n_4$ increases, the cost of solving the FOM increases with $\mathcal{O}(n^4)$ while the cost of \texttt{TT-CUR-DEIM} scales linearly.   As an example, for the 4D-AR case with  $n=81$, FOM takes about 10.5 hours to complete 4000 time-steps while \texttt{TT-CUR-DEIM} takes only 1.8 minutes. Also, for S-ADR problem with \(n_2 = n_3 = 50\), the FOM requires about 21.9 hours to execute, whereas \texttt{TT-CUR-DEIM} takes approximately 16.7 minutes for 10,000 time steps.

The required memory for storing the RHS of the \texttt{TT-CUR-DEIM} model relative to the RHS of the FOM is also compared in Figure \ref{fig:Barchart} for 4D-AR and S-ADR examples.

\section{Conclusion}
We present a novel algorithm for the time integration of tensor differential equations on low-rank manifolds. This algorithm is applicable to both tensor train and Tucker tensor low-rank approximations. It advances the low-rank tensor by applying cross approximation algorithms to the time-discrete full-order model. These algorithms, which are based on the Discrete Empirical Interpolation Method (DEIM), sample the time-discrete full-order model at strategically selected entries. For the Tucker tensor low-rank approximation, we utilize \texttt{DEIM-FS} that was recently developed \cite{GB24}. For the tensor train low-rank approximation, we introduce \texttt{TT-CUR-DEIM}, a tensor train cross approximation based on DEIM.

The algorithm offers multiple advantages: (i) Solving nonlinear TDEs on tensor train low-rank manifolds requires evaluating the full-order model at $\mathcal{O}(r^2dn)$ entries, which is cost-optimal, i.e., the minimum number of function evaluation to construct a rank-$r$ tensor train approximation. For a rank-$r$ Tucker tensor decomposition, \texttt{DEIM-FS} requires $\mathcal{O}(r^d+ ndr^{d-1})$ entries, which is significantly smaller than $n^d$. The algorithm is agnostic to the type of nonlinearity.  (ii) The time-integration scheme does not utilize the inversion of auto-correlation matrices and it is robust in the presence of small or zero singular values. (iii) We demonstrate that the presented algorithm inherits the temporal accuracy of the time-discrete full-order model. In particular, we develop high-order Runge-Kutta time integration schemes. (iii) The algorithm is simple to implement as it does not utilize tangent space projection and it does require a term-by-term treatment of the right-hand-side terms. 

 We demonstrate the performance of the algorithm for several test cases including a 100-dimensional nonlinear TDE. The Tucker tensor and tensor train low-rank decomposition demonstrate the potential applicability of the presented algorithm to tensor networks. 

\enlargethispage{20pt}

\section*{Acknowledgments}
We thank Professor Peyman Givi for many stimulating discussions that led to a number of improvements in this paper. This work is sponsored by the National Science Foundation (NSF), USA under Grant CBET-2152803 and by the Air Force Office of Scientific Research award no. FA9550-21-1-0247.

\vskip2pc

\appendix

\section{DEIM algorithm}
The DEIM pseudocode is presented via Algorithm \ref{alg:DEIM}. This algorithm is adopted from \cite{SCNonlinear}.

\begin{algorithm}
\begingroup
\fontsize{9pt}{9pt}\selectfont
\SetAlgoLined
\KwIn{$\mathbf{U}_{p}=\left[\begin{array}{llll}\mathbf{u}_{1} & \mathbf{u}_{2} & \cdots & \mathbf{u}_{p}\end{array}\right]$}
\KwOut{$\bm I_{p}$}
$\left[\rho, \bm I_{1}\right]=\max \left|\mathbf{u}_{1}\right|$ \hspace{15mm} $\rhd$ choose the first index\;
$\mathbf{P}_{1}=\left[\mathbf{e}_{\bm I_{1}}\right]$ \hspace{24mm} $\rhd$ construct first measurement matrix\;
\For{$i=2$ \KwTo $p$}{
$\mathbf{P}_{i}^{T} \mathbf{U}_{i} \mathbf{c}_{i}=\mathbf{P}_{i}^{T} \mathbf{u}_{i+1}$ \hspace{7.5mm}  $\rhd$ calculate $c_{i}$\;
$\mathbf{R}_{i+1}=\mathbf{u}_{i+1}-\mathbf{U}_{i} \mathbf{c}_{i}$ \hspace{7mm}  $\rhd$ compute residual\;
$\left[\rho, \bm I_{i}\right]=\max \left|\mathbf{R}_{i+1}\right|$ \hspace{7mm} $\rhd$ find index of maximum residual\;
$\mathbf{P}_{i+1}=\left[\begin{array}{ll}\mathbf{P}_{i} & \mathbf{e}_{\bm I_{i}}\end{array}\right]$ \hspace{6mm} $\rhd$ add new column to measurement matrix\;}
\caption{\texttt{DEIM} Algorithm \cite{SCNonlinear}}
\label{alg:DEIM}

\endgroup
\end{algorithm}

\bibliographystyle{ieeetr}
\bibliography{References,Hessam}

\end{document}